\documentclass[11pt,a4]{amsart}
\usepackage{amsthm,amssymb,amsfonts,amsmath,color,bm}
\usepackage{graphicx,enumerate}
\usepackage{comment,epsfig}
\usepackage{multirow}
\usepackage{mathrsfs}
\usepackage{caption}
\usepackage{amsaddr}
\usepackage{type1cm}

\usepackage[bottom]{footmisc}
\newcommand{\bfd}{{\boldsymbol d}}
\newcommand{\bfxi}{{\boldsymbol \xi}}

\renewcommand{\a}{\alpha}
\newcommand{\red}{\color{red}}
\newcommand{\blue}{\color{blue}}
\newcommand{\green}{\color{carageen}}

\newcommand{\clt}{central limit theorem}

\newcommand{\garch}{{\rm GARCH}$(1,1)$}

\newcommand{\ex}{{\rm e}\,}

\newcommand{\asy}{asymptotic}

\newcommand{\ts}{time series}
\newcommand{\tsa}{\ts\ analysis}

\newtheorem{lemma}{Lemma}[section]

\newtheorem{theorem}[lemma]{Theorem}

\renewcommand{\P}{{\mathbb P}}
\newtheorem{proposition}[lemma]{Proposition}
\newtheorem{definition}[lemma]{Definition}
\newtheorem{corollary}[lemma]{Corollary}
\newtheorem{example}[lemma]{Example}
\newtheorem{exercise}[lemma]{Exercise}
\newtheorem{remark}[lemma]{Remark}

\newtheorem{tab}[lemma]{Table}

\newcommand{\MC}{Markov chain}

\newcommand{\bfQ}{{\bf Q}}

\newcommand{\bfu}{{\bf u}}

\newcommand{\bfR}{{\bf R}}

\newcommand{\bth}{\begin{theorem}}
\newcommand{\ethe}{\end{theorem}}
\newcommand{\sv}{stochastic volatility}
\newcommand{\bre}{\begin{remark}\em }
\newcommand{\ere}{\end{remark}}

\newcommand{\ble}{\begin{lemma}}
\newcommand{\ele}{\end{lemma}}
\newcommand{\sre}{stochastic recurrence equation}
\newcommand{\pp}{point process}
\newcommand{\bde}{\begin{definition}}
\newcommand{\ede}{\end{definition}}
\newcommand{\bco}{\begin{corollary}}
\newcommand{\eco}{\end{corollary}}
\newcommand{\bpr}{\begin{proposition}}
\newcommand{\epr}{\end{proposition}}

\newcommand{\bexer}{\begin{exercise}}
\newcommand{\eexer}{\end{exercise}}

\newcommand{\bexam}{\begin{example}}
\newcommand{\eexam}{\end{example}}

\newcommand{\efi}{\end{fig}}

\newcommand{\btab}{\begin{tab}}
\newcommand{\etab}{\end{tab}}

\newcommand{\lhs}{left-hand side}

\newcommand{\rv}{random variable}

\newcommand{\sign}{{\rm sign}}

\newcommand{\var}{{\rm var}}

\newcommand{\as}{{\rm a.s.}}

\newcommand{\bfTh}{\mbox{\boldmath$\Theta$}}

\newcommand{\rhs}{right-hand side}
\newcommand{\df}{distribution function}

\newcommand{\dint}{\displaystyle\int}

\newcommand{\beao}{\begin{eqnarray*}}
\newcommand{\eeao}{\end{eqnarray*}\noindent}

\newcommand{\beam}{\begin{eqnarray}}
\newcommand{\eeam}{\end{eqnarray}\noindent}

\newcommand{\beqq}{\begin{equation}}
\newcommand{\eeqq}{\end{equation}\noindent}

\newcommand{\bce}{\begin{center}}
\newcommand{\ece}{\end{center}}

\newcommand{\barr}{\begin{array}}
\newcommand{\earr}{\end{array}}

\newcommand{\std}{\stackrel{d}{\rightarrow}}
\newcommand{\stas}{\stackrel{\rm a.s.}{\rightarrow}}

\newcommand{\stw}{\stackrel{w}{\rightarrow}}

\newcommand{\eqd}{\stackrel{d}{=}}

\newcommand{\vague}{\stackrel{\lower0.2ex\hbox{$\scriptscriptstyle
                    \it{v} $}}{\rightarrow}}
\newcommand{\weak}{\stackrel{\lower0.2ex\hbox{$\scriptscriptstyle
                    \it{w} $}}{\rightarrow}}
\newcommand{\what}{\stackrel{\lower0.2ex\hbox{$\scriptscriptstyle
                    \it{\hat{w}} $}}{\rightarrow}}

\newcommand{\bdis}{\begin{displaymath}}
\newcommand{\edis}{\end{displaymath}\noindent}

\newcommand{\R}{\mathbb{R}}

\newcommand{\nto}{n\to\infty}

\newcommand{\xto}{x\to\infty}

\newcommand{\ov}{\overline}
\newcommand{\wt}{\widetilde}

\newcommand{\vep}{\varepsilon}

\newcommand{\la}{\lambda}

\newcommand{\regvary}{regularly varying}
\newcommand{\slvary}{slowly varying}
\newcommand{\regvar}{regular variation}

\newcommand{\bbr}{{\mathbb R}}

\newcommand{\bbz}{{\mathbb Z}}
\newcommand{\Z}{{\mathbb Z}}

\newcommand{\con}{convergence}

\newcommand{\evt}{extreme value theory}

\newcommand{\st}{such that}
\newcommand{\fif}{if and only if}
\newcommand{\wrt}{with respect to}
\newcommand{\chf}{characteristic function}
\newcommand{\chdf}{characteristic-\ds\ \fct }
\newcommand{\fct}{function}

\newcommand{\ds}{distribution}

\newcommand{\rep}{representation}

\newcommand{\seq}{sequence}

\newcommand{\pro}{probabilit}

\newcommand{\ms}{measure}
\newcommand{\mgf}{moment generating function}

\newcommand{\bfx}{{\bf x}}
\newcommand{\bfX}{{\bf X}}

\newcommand{\bfS}{{\bf S}}
\newcommand{\bfW}{{\bf W}}

\newcommand{\bfq}{{\bf q}}

\newcommand{\E }{{\mathbb E}}
\renewcommand{\P }{{\mathbb P}}

\newcommand{\1}{{\mathbf 1}}

\allowdisplaybreaks

\begin{document}
\today
\title{Moments for self-normalized partial sums}
\author[Matsui, Mikosch, Wintenberger]{Muneya Matsui, Thomas Mikosch, 
Olivier Wintenberger}
\address{M. Matsui\\ 
Department of Business Administration, 
Nanzan University\\
18 Yamazato-cho Showa-ku Nagoya, 466-8673, Japan. \\
\small{mmuneya@gmail.com}
}
\address{T. Mikosch\\
Department  of Mathematics,
University of Copenhagen \\
Universitetsparken 5,
DK-2100 Copenhagen,
Denmark. \\
\small{mikosch@math.ku.dk}
}
\address{O. Wintenberger\\
Laboratoire de Probabilit\'es, Statistique et Mod\'elisation\\
Sorbonne Universit\'e, UPMC Université Paris 06, F-75005, Paris, France\\
and\\
Institut CNRS Pauli, Vienna University. \\
\small{olivier.wintenberger@sorbonne-universite.fr}
}

\begin{abstract}
We consider a \regvary\ stationary 
\seq\ of random variables $(X_t)$ with tail index $\a<2$. 
For these \seq s we  study the joint \con\ of sums, $\ell^p$-type moduli and 
maxima. We focus on ratio statistics, 
including the studentized sums and sums normalized by the corresponding maxima, and study the existence
of moments for the limit ratios. We consider particular examples of processes 
$(X_t)$  whose limit ratios possess all moments. 
But, in contrast to the latter situation, there also exist \seq s $(X_t)$
where certain moments of the limit ratio are infinite. This phenomenon results from extremal clusters in the \seq .
\end{abstract}


\maketitle
\par\noindent
{\em MSC2020 subject classifications:}{\rm Primary 60F05; Secondary 60E10 60G70 62E20}


\section{Introduction}\setcounter{equation}{0}
Standardization
is an essential step when pre-processing data, in particular when
dealing with machine learning models and statistical analyses. This standardization
ensures that all features of the model 
contribute equally to the analysis by transforming them to common scale, 
thereby enhancing both the performance and interpretability of the results. 
Additionally, standardization can help to reduce the impact of extreme values, 
making their influence less pronounced and leading to more robust model 
performance. In this paper, we focus on examining the impact of 
dependence among extreme values on the standardization process. 
We consider the moment properties of the limit under 
self-normalization for partial sums that 
are attracted toward $\a$-stable distributions for $\a<2$; for an encyclopedic study
in the case of Gaussian limits
we refer to the monograph by de la Pe\~{n}a et al. \cite{pena:lai:shao:2009}.
\subsection{Infinite-variance \clt s}
We consider an $\bbr^d$-valued (strictly) stationary  
\seq\ $(\bfX_t)$ with generic element $\bfX$ 
and the corresponding partial sum process 
\beao
\bfS_n=\bfX_1+\cdots +\bfX_n\,,\qquad n\ge 1\,. 
\eeao
Assuming finite variance of $\bfX$ and suitable dependence conditions on $(\bfX_t)$,
central limit theory
for $(\bfS_n)$ is well studied. The infinite-variance case 
has attracted less attention, the papers by 
Jakubowski \cite{jakubowski:1993,jakubowski:1997},
Davis and Hsing \cite{davis:hsing:1995} being exceptions. 
In addition to infinite second moments one needs
to specify the joint \asy\ tail behavior in the \seq\  $(\bfX_t)$. 
A natural condition in this context is {\em serial \regvar }.
Following Basrak and Segers \cite{basrak:segers:2009}, $(\bfX_t)$ is \regvary\
with index $\a>0$
if there exist a Pareto$(\a)$-distributed \rv\ $Y_\a$, i.e., $\P(Y_\a>x)=x^{-\a}$, 
$x>0$,  
independent of an 
$\bbr^d$-valued \seq\ $(\bfTh_t)$ (the {\em spectral tail process}) \st\ for every $h\ge 0$, 
\beam\label{def:regvar}
\P\big(x^{-1}(\bfX_{-h},\ldots,\bfX_h)\in\cdot \,\big|\,|\bfX_0|>x\big)
\stw \P\big(Y_\a\,(\bfTh_{-h},\ldots,\bfTh_h)\in \cdot\big)\,,\quad \xto\,.
\eeam
This means that for every $h\ge 0$, $(\bfX_{-h},\ldots,\bfX_h)$ satisfies the 
condition of multivariate \regvar\ in $\bbr^{d(h+1)}$ with index $\a$; 
cf. Resnick \cite{resnick:1987,resnick:2007}. For an iid 
\seq\ $(\bfX_t)$ Rva\v ceva \cite{rvaceva:1962} proved that \eqref{def:regvar}
with $h=0$ is necessary and sufficient for the $\a$-stable \clt , $\a\in (0,2)$,
\beam\label{eq:stabclt}
a_n^{-1}\big(\bfS_n-\bfd_n)\std \bfxi_\a\,,\qquad \nto\,,
\eeam  
where $\bfxi_\a$ is an $\bbr^d$-valued $\a$-stable random vector (see
Feller \cite{feller:1971}, Samorodnitsky and Taqqu \cite{samorodnitsky:taqqu:1994}
for the definition of an $\a$-stable random vector), the normalizing 
\seq\ $(a_n)$ satisfies $n\,\P(|\bfX|>a_n)\to 1$ as $\nto$, and $(\bfd_n)$ are 
suitable centering constants. 
In the iid case, 
\eqref{def:regvar} for $h=0$ and for any $h\ge 1$ are equivalent.
Under the serial \regvar\ condition \eqref{def:regvar}
extensions of \eqref{eq:stabclt} to stationary \seq s $(\bfX_t)$ were 
given in 
\cite{jakubowski:1993,jakubowski:1997,davis:hsing:1995,basrak:segers:2009,bartkiewicz:jakubowski:mikosch:wintenberger:2011}, where the  limit vector $\bfxi_\a$ can be  expressed in terms of $\a$ and $(\bfTh_t)$. Recently, various books on  \regvary\
\ts\ models have been published: Buraczewski et al. \cite{buraczewski:damek:mikosch:2016}, Kulik and Soulier \cite{kulik:soulier:2020}, Mikosch and 
Wintenberger \cite{mikosch:wintenberger:2024}. Among others, 
these texts  discuss conditions for the \clt\ \eqref{eq:stabclt}, \evt , \pp\ \con ,
clustering phenomena, and they demonstrate the  serial \regvar\ property 
for a variety of
stationary \ts\ models, including the classical linear processes with iid \regvary\ noise, affine \sre s and the GARCH model, max-moving averages, max-stable processes, and the Lindley process. 
\subsection{Self-normalization of partial sums}
Chapter 10 in \cite{mikosch:wintenberger:2024} and the recent paper
Matsui et al. \cite{matsui:mikosch:wintenberger:2024} have been devoted to
limit theory for $(\bfS_n)$ under  {\em self-normalization}. Self-normalization 
is a classical topic in \pro y theory and statistics; see e.g. the textbook
by de la Pe\~{n}a et al. \cite{pena:lai:shao:2009}. In the infinite-variance 
case there is a particular need for self-normalizations because the normalizing 
\seq\ $(a_n)$ in \eqref{eq:stabclt} depends on the tail of $|\bfX|$ which is
typically unknown. Therefore it is desirable to replace $a_n$ by a {\em known} random
quantity of the same magnitude. For $\a\in (0,2)$ and under
additional assumptions 
the {\em  maximum} $\gamma_{n,\infty}=\max_{t=1,\ldots,n}|\bfX_t|$, $n\ge 1$, 
and the {\em $\ell^p$-moduli} for $p>\a$,
\beao
\gamma_{n,p}=\Big(\sum_{t=1}^n |\bfX_t|^p\Big)^{1/p}\,,\qquad n\ge 1\,,
\eeao
have the same \asy\ order as $\bfS_n$, i.e.,
\beao
a_n^{-1}\big(\bfS_n-\bfd_n,\gamma_{n,\infty},\gamma_{n,p}\big)\std \big(\bfxi_\a,\zeta_{\a,\infty},\zeta_{\a,p}\big)\,, \qquad \nto\,,
\eeao 
where $\bfxi_\a$ is $\a$-stable, $\zeta_{\a,\infty}$ is Fr\'echet-distributed with shape
parameter $\a$,
and $\zeta_{\a,p}^p$ is a positive $\a/p$-stable \rv .
Hence as $\nto$,
\beam\label{eq:selfnorm1}
\dfrac{\bfS_n-\bfd_n}{\gamma_{n,\infty}}&\std& \dfrac{\bfxi_\a}{\zeta_{\a,\infty}}=:\bfR_{\a,\infty}\,,\\  
\dfrac{\bfS_n-\bfd_n}{\gamma_{n,p}}&\std& \dfrac{\bfxi_\a}{\zeta_{\a,p}} =:\bfR_{\a,p}\,,\quad p>\a\,,\label{eq:selfnorm2}
\eeam
but also $\gamma_{n,p}/\gamma_{n,\infty}\std \zeta_{\a,p}/\zeta_{\a,\infty}$. In particular, the case $p=2$
corresponds to the classical studentization. 
\par
The limit relations \eqref{eq:selfnorm1}, \eqref{eq:selfnorm2} solve the problem of replacement of the 
unknown $a_n$ by a known self-normalization. But they also create another difficulty: the limit ratio \ds s are not familiar; only little is known about their
properties. One goal of self-normalizing procedures is 
that $q$th moments of the limit ratios for $q>\a$ might exist. Since 
the stable vector $\bfxi_\a$ has only moments of order $q<\a$, higher moments
of the limit ratios indicate that their \ds s are less spread out than those
of $\bfxi_\a$. An objective of our paper is to show that this goal can often be
achieved. But we will also prove that, in contrast to an iid \seq\ $(\bfX_t)$,
infinite power-moments of the limit ratios are possible due to extremal clustering
effects in the stationary \seq .   
\par
We also intend to give more precise \ds al information about the limit ratios 
than the existence of certain moments. For this reason we start by recalling
some classical results in the iid real-valued case.

\subsection{The iid benchmark}\label{subsec:iid} For comparison 
with the dependent case we collect some results on self-normalizing quantities
under the assumption that $(X_t)$ is an iid real-valued \seq\ and \regvary\ with index $\a\in (0,2)$. In particular, we call a generic element $X$ {\em \regvary\ with index $\a>0$} if for some \slvary\ \fct\ $L$,
\beao
 \P(|X|>x)=L(x)x^{-\a}\,,x>0\,,\mbox{ and the limits } \lim_{\xto}\dfrac{\P(\pm X>x)}{\P(|X|>x)}=p_{\pm} \mbox{ exist.}
\eeao 
The {\em tail-balance coefficients} $p_{\pm}$ define the 
\ds\ of $\Theta_0$: $\P(\Theta_0=\pm 1)=p_{\pm}$.
\par
Chistyakov and G\"otze \cite{chistyakov:goetze:2004} considered the case $p=2$.
They proved that $S_n/\gamma_{n,2}\std Z$ for some genuine \rv\ $Z$ implies that $Z$ is Gaussian or 
$X$ is \regvary\ with index $\a\in (0,2)$. Moreover, in the latter case and for $\a\in (1,2)$, $\E[X]=0$ holds necessarily.
\par
Logan et al. \cite{logan:mallows:rice:shepp:1973} investigated
the limits \eqref{eq:selfnorm2} for $p>\a$. In particular, they found 
the densities of the limit ratios $R_{\a,p}$ for $p=1$ and $p=2$
and showed that they have all moments. For $p=2$, $\a\in (1,2)$ the densities are \asy ally equivalent to
Gaussian densities. If $p=1$,   $\a<1$ then $|S_n|/\gamma_{n,1}\le 1$ a.s., and if
$p_+p_->0$ then  $ R_{\a,1}$ has a continuous density on $(-1,1)$
which has singularities at~$\pm 1$. 
For iid \rv s $(X_t)$ that are also positive for $\a\in (0,1)$ and centered positive for $\a\in (1,2)$,
the limit \ds\ of $ R_{\a,\infty}$ was derived by Darling \cite{darling:1952}. 
The general case was considered in Matsui et al. \cite{matsui:mikosch:wintenberger:2024}.

\ble\label{lem:1} Assume that $X$ is \regvary\ with index 
$\a\in (0,2)\backslash\{1\}$ and $\E[X]=0$ for $\a\in (1,2)$. Then
 $S_n/\gamma_{n,\infty}\std R_{\a,\infty}$, $\nto$, and the limit ratio has a finite \mgf\
in some neighborhood of the origin.  
\ele
\begin{proof} Corollary 4.1 of \cite{matsui:mikosch:wintenberger:2024}
yields
the \chf 
\beam\label{eq:mgffinite}
\varphi_{ R_{\a,\infty}}(u)&=& \E\big[\ex^{i\,u\,R_{\a,\infty} }\big]\nonumber\\
&=&\dfrac{p_+\,\ex^{iu}+p_-\,\ex^{-iu}}{
1+ iu \frac \a {\a-1}\,(p_+-p_-)\,\1_{(1,2)}(\a)+ \int_{[-1,1]}\big( 1+iyu\,\1_{(1,2)}(\a)-
\ex^{iyu}\big)\,\mu_\a(dy) }\,,
\eeam
where $\mu_\a(dy) =\big(p_+\,\1(y>0)+p_-\,1(y<0)\big)\alpha/|y|^{\alpha+1} dy\,.$
We will show that the \mgf\ $\varphi_{R_{\a,\infty}}(-i\,u)$ is finite 
for $u\in (-\vep,\vep)$ for some small $\vep>0$. We focus on the denominator
on the \rhs\ of \eqref{eq:mgffinite}. For $\a\in (1,2)$, 
$1+ u \frac \a {\a-1}\,(p_+-p_-)$ is arbitrarily close to 1 for small
values of $|u|$. On the other hand, a Taylor expansion argument shows that
for small $|u|$ and some constant $c>0$,
\beao
&& \Big| \int_{[-1,1]}\big( 1+iyu\,\1_{(1,2)}(\a)-
\ex^{iyu}\big)\,\mu_\a(dy)\Big| \\
&& \le 
c\,|u|\,\int_{[-1,1]} \big(\1_{(1,2)}(\a)\,
|y|^{1-\a}+\1_{(0,1)}(\a)\,|y|^{-\a}\big)\,dy\,.
\eeao
The \rhs\ is arbitrarily small for small values of $|u|$. This finishes the
proof. 
\end{proof}
\par
The results in the iid case show that self-normalizations 
can be beneficial for statistical analyses. Indeed, for $\a\in (0,2)$, $p>\a$, 
the tails of $R_{\a,p}$ are significantly lighter than the classical
$\a$-stable limits $\xi_\a$ of $(a_n^{-1}S_n)$ that exhibit power-law tails
of order $\a$. However, this (positive) result does not apply to all extremal dependence structures. 
For example, Anderson and Turkman \cite{anderson:turkman:1995} presented a real-valued \regvary\ stationary process \st\ the components of the \ds al limit of 
$a_n^{-1}(S_n\,,\max_{1\le t\le n}X_t)$ are independent, implying infinite moments of order $\a$ for the limit ratio. Recently, Matsui et al.  \cite{matsui:mikosch:wintenberger:2024}  proved the integrability of the limit ratio of  $(S_n/\gamma_{n,\infty})$ under general extremal dependence 
conditions, even for $\a<1$. Thus the standardization of the sums $S_n$ by $\gamma_{n,\infty}$ is beneficial, and even more so if one chooses
$\gamma_{n,p}$ for finite $p>\a$, such as for the classical studentization with $p=2$. 
\subsection{Goals of this paper} In what follows, we will study properties of the \asy\ 
\ds\ of the $\bbr^d$-valued self-normalized quantity
$\bfS_n/\gamma_{n,p}$ for $\a\in (0,2)$ and $p>\a$. In particular, we are interested in the existence of moments of these quantities.
In most parts of the paper (for example, when looking at examples) we
will focus on univariate \ts , i.e., $d=1$. Then we will not use  boldface symbols. 
\par
In Section~\ref{sec:2} we introduce {\em mixing} and {\em anti-clustering conditions} for \regvary\ \ts . From Matsui et al. \cite{matsui:mikosch:wintenberger:2024} we recall 
Theorem~\ref{pr:prophybridch} about the joint \con\ of maxima, sums and $\ell^p$-moduli. This result provides the joint \chf-\df-Laplace transform
of these quantities; it is the main tool for deriving quantitative descriptions of the limit ratios of self-normalized expressions. Another useful tool is
\pp\ \rep s of the \ds al limits for sums and $\ell^p$-moduli. These are in the spirit of infinite series \rep s of $\a$-stable  random vectors 
(as given in Samorodnitsky and Taqqu \cite{samorodnitsky:taqqu:1994}) and of \rep s of extremes (as given in Davis and Hsing \cite{davis:hsing:1995}, Basrak and Segers 
\cite{basrak:segers:2009}, Kulik and Soulier \cite{kulik:soulier:2020}).
\par
Our main results are presented in Section~\ref{sec:mom}, in particular in Theorem~\ref{thm:momentsratio}.
The latter provides necessary and/or sufficient conditions for the existence of integer-moments of the limit ratios $\bfR_{\a,p}$, $p>\a$.
 The \ds s of $ \bfR_{\a,p}$  are mixtures of
the limit vector in the iid case and quantities that depend only on the spectral
tail process. In Section~\ref{sec:exam} 
we apply these results to various univariate \regvary\ stationary \ts\ models.  
In certain cases (such as linear processes, solutions to affine \sre s, \sv\  and GARCH processes) the \ds\ of $ R_{\a,p}$
is "well-behaved" in the sense that all moments are finite, as for the 
iid benchmark.
\par
The \ds\  of $R_{\a,p}$  is the same as in the iid case
(or scaled version of those) \fif\ $(\Theta_t)_{t\ne 0}$ is deterministic. These 
include the cases of \regvary\ \sv\ models (when $\Theta_t=0$ a.s. for $t\ne 0$; 
see Section~5.3 in \cite{mikosch:wintenberger:2024})
and \regvary\ linear processes (when $\Theta_t$, $t\ne 0$, are deterministic; see
Section~5.2 in \cite{mikosch:wintenberger:2024}). For GARCH processes and solutions to affine stochastic recurrence equations the limit ratios have distinct \ds s.
But even worse, in Section \ref{sec:exam} we also construct a \regvary\ stationary \ts\ $(X_t)$ that satisfies all conditions of Theorem~\ref{thm:momentsratio} but certain moments of $ R_{\a,p}$ are infinite.
This lack of moments for $R_{\a,p}$ may occur for non-negative time series 
$(X_t)$ where the appearance of a large block $\sum_{t=s}^r X_t$ does not necessarily coincide with an extreme $X_t$ at a certain time $s\le t\le r$. 
Such models are interesting in themselves -- as processes with a positive extremal index but with potentially rather long periods of extremal dependence which impact the 
moment properties of self-normalizations and where standardization fails.
\par
Finally, Section~\ref{sec:even:moments} is entirely devoted to the derivation of  sufficient conditions for the existence
of  even moments 
of the limit ratio $R_{\a,2}$, i.e., to the limit of studentized sums, and to the calculation of these moments. 
Appendix A supplements this section by calculating and illustrating the density of $R_{\a,2}$ in the case of a \regvary\ iid 
\seq\ with index $\a\in (1,2)$. Here we closely follow the classical paper Logan et al. \cite{logan:mallows:rice:shepp:1973}
that contains the surprising result that the density of the limit ratio has tails that are equivalent  to the Gaussian density. This nice property extends from the iid case  to \sv\  and linear models.

\section{Preliminaries}\label{sec:2}\setcounter{equation}{0}
Throughout we consider an $\bbr^d$-valued \regvary\ stationary process $(\bfX_t)$ 
with index $\alpha\in  (0,2)\backslash \{1\}$. If $\a\in (1,2)$ we also assume 
that $\E[\bfX]=0$, and if we restrict ourselves to the case $d=1$ we will not use boldface symbols.
We choose the normalizing constants
$(a_n)$ \st\
$n\,\P(|\bfX|> a_n)\to 1$ as $\nto$. 
\par
We introduce the  {\em spectral cluster process} 
$\bfQ=\bfTh/\|\bfTh\|_\a$, where we define $\|\cdot\|_p$, $p\in[\a,\infty]$, such that $\|\bfx\|_p^p=\sum_{t\in\Z}|\bfx_t|^p$ for $\bfx \in \ell^\a(\bbr^d)$.
The a.s. \con\ of the series $\|\bfTh\|_\a^\a$ is always guaranteed under the 
conditions  introduced below; see Janssen \cite{janssen:2019}, cf. Section 6.2 in 
\cite{mikosch:wintenberger:2024}.
\subsection{Mixing conditions} 
We will use the following {\bf mixing condition}:
for some integer \seq s $r_n\to \infty$, $k_n=[n/r_n]\to\infty$,  
\beam\label{eq:wdepms2}
\Psi_{n,p}(\bfu,x,\la)&:=&\E\Big[\exp\big(i\,a_n^{-1}\bfu^\top\bfS_{n}-a_n^{-p}\la \gamma_{n,p}^p\big)\,\1\big(a_n^{-1}{\red \gamma_{n,\infty}}
\le x\big)\Big]\nonumber\\
 &=&
\Big(\E\Big[\exp\big(i\,a_n^{-1}\bfu^\top\bfS_{r_n}-a_n^{-p}\la \gamma_{r_n,p}^p\big)\,\1\big(a_n^{-1}{\red \gamma_{r_n,\infty}}
\le x\big)\Big]\Big)^{k_n}+o(1),\nonumber\\&&\qquad \nto,\qquad (\bfu,x,\la)\in\bbr^d\times \R_+^2\,.
\eeam
This mild condition ensures the \asy\ independence of $k_n$ maxima, sums, and 
$\ell^p$-moduli over disjoint blocks of length $r_n$. It
follows from strong mixing properties of $(\bfX_n)$ or by coupling arguments; see for example Rio \cite{rio:2017}, the discussion of this condition  
in Chapters 9 and 10
in \cite{mikosch:wintenberger:2024}, or in Matsui et al. \cite{matsui:mikosch:wintenberger:2024}.
\subsection{Anti-clustering condition}
We will use the following {\bf anti-clustering condition}: for some $r_n\to\infty$ \st\ $k_n=[n/r_n]\to\infty$,
\beam\label{cond:acac}
\lim_{l\to\infty}\limsup_{\nto}  n \,\sum_{j=l}^{r_n}\E\big[(|a_n^{-1}\bfX_{j} |\wedge  1)\,(|  a_n^{-1}\bfX_{0}|\wedge
1)\big]=0\,.
\eeam
This condition ensures that joint high-level exceedances of $|\bfX_0/a_n|$ and $|\bfX_j/a_n|$ for $j\in [l,r_n]$ and sufficiently large $l$ are unlikely.
This condition is satisfied for a large variety of \regvary\ \ts\ models; see
Chapter 9 and 10 in \cite{mikosch:wintenberger:2024} and Matsui et al.
\cite{matsui:mikosch:wintenberger:2024}. 

\subsection{Joint \con\ of sums, maxima, and $\ell^p$-moduli}
We start with a limit result about the joint \con\ of $a_n^{-1}(\bfS_n,\gamma_{n,\infty},\gamma_{n,p})$ as $\nto$. It was proved in Matsui et al. \cite{matsui:mikosch:wintenberger:2024}.
\bth\label{pr:prophybridch}
Assume the conditions \eqref{eq:wdepms2} and \eqref{cond:acac}. 
Then, with the notation introduced above,  for $\a<p$,
\beao
a_n^{-1} (\bfS_n, \gamma_{n,\infty},\gamma_{n,p})\std (\bfxi_\a,\zeta_{\a,\infty},\zeta_{\a,p})\,,\qquad \nto\,,
\eeao
where the joint limit \ds\  
is described by the joint hybrid \chdf --Laplace transform
of  $(\bfxi_\a,\zeta_{\a,\infty},\zeta_{\a,p}^p)$ given by
\beam\label{eq:xxx}\Psi_\bfX(\bfu,x,\la)&:=&
\E\big[\ex^{i \,\bfu^\top\,\bfxi_\a}\,\1(\zeta_{\a,\infty}\le x)\,\ex^{-\la\,\zeta_{\a,p}^p}\big]\nonumber\\
&=&\exp\Big(\dint_0^\infty\E\Big[\ex^{i\,y\, \bfu^\top\sum_{t=-\infty}^\infty  \bfQ_t-y^p\la\sum_{t=-\infty}^\infty  |\bfQ_t|^p}\, 
\1\Big(y\,\max_{t\in \Z}|\bfQ_t|\le x \Big)\nonumber\\
&&\nonumber-1-
i\,y\,\bfu^\top \sum_{t\in\Z} \bfQ_t\,\1_{(1,2)}(\a)\Big]\,d(-y^{-\a})\Big)\,,\qquad (\bfu,x,\la)\in\bbr^d\times \R_+^2\,.
\eeam
\ethe
\bre\label{rem:moments}
The limit \ds\ involves the infinite series $\bfW=\sum_{t\in\bbz}\bfQ_t$. It converges 
absolutely a.s. and $\E[(\sum_{t\in\bbz}|\bfQ_t|)^\a]<\infty$. For $\a\in (0,1)$ this is a trivial con\seq\ of the fact that
$\|\bfTh\|_\a<\infty$ and  $\|\bfQ\|_\a=1$ a.s. by definition of $\bfQ$; see Janssen \cite{janssen:2019}. 
For $\a\in (1,2)$ these properties are con\seq  s of the anti-clustering condition \eqref{cond:acac}; see Lemma 9.1.7 in 
\cite{mikosch:wintenberger:2024}.
\ere
\bre
For $\la=0,x=\infty$ we recognize the \chf\ of an $\a$-stable vector $\bfxi_\a$:
\beao
\varphi_{\bfxi_\a}(\bfu)&=&\E\big[\ex^{i \,\bfu^\top\,\bfxi_\a}\big]\\&=&
\exp\Big(\dint_0^\infty\E\big[\ex^{i\,y\, \bfu^\top\bfW}
-1
-
i\,y\,\bfu^\top \bfW\,\1_{(1,2)}(\a)\big]\,d(-y^{-\a})\Big)\,,\qquad \bfu\in\R^d\,.
\eeao
For $\la=0,\bfu=\bf0$ we obtain the \df\ of $\zeta_{\a,\infty}$:
\beao
\P(\zeta_{\a,\infty}\le x)
&=& \red \exp\Big(-\int_0^\infty \P\Big(y\,\max_{t\in \Z}|\bfQ_t|>x \Big)\,d(-y^{-\a})\Big)=\ex^{-x^{-\a} \,\theta_{|\bfX|}}\,,\qquad x>0\,,
\eeao
where $\theta_{|\bfX|}=\E\big[\max_{t\in \Z}|\bfQ_t|^\a \big]$ is the
extremal index of the \seq\ $(|\bfX_t|)$; see Davis and Hsing \cite{davis:hsing:1995}, Basrak and Segers \cite{basrak:segers:2009}. Under our conditions, 
$\theta_{|\bfX|} \in (0,1]$; see \cite{basrak:segers:2009}; cf. Corollary~6.2.2
in \cite{mikosch:wintenberger:2024}. We notice that $\zeta_{\a,\infty}$ has a Fr\'echet
\ds .\\[1mm]
For $\bfu=\bf0$, $x=\infty$ we obtain the Laplace transform of a positive $\a/p$-stable \rv\ $\zeta_{\a,p}^p$:
\beao
\E\big[\ex^{-\la\,\zeta_{\a,p}^p}\big]
&=&\exp\Big(\dint_0^\infty\E\Big[\ex^{-y^p\la\sum_{t=-\infty}^\infty  |\bfQ_t|^p}-1\Big]\,d(-y^{-\a})\Big)\\
&=&\exp\big( -\E\big[Y_p^\a]\E\big[\|\bfQ\|_p^\a\big]\la^{\a/p} \big)\,,
\eeao
where $Y_p$ is Pareto$(p)$-distributed; see Matsui et al. \cite{matsui:mikosch:wintenberger:2024}.
\ere
By virtue of Theorem~\ref{pr:prophybridch} we have the joint \con\ of the 
normalized vectors $a_n^{-1}(\bfS_n,\gamma_{n,\infty},\gamma_{n,p})$ for $p>\a$. Therefore the desired 
limit ratios in \eqref{eq:selfnorm1} and \eqref{eq:selfnorm2} exist.
\subsection{Point process \rep\ of the limits}
In what follows, we will express the joint limit in Theorem~\ref{pr:prophybridch} in terms of Poisson
points. We consider an enumeration of the points $\Gamma_1< \Gamma_2<\cdots$
of a unit rate Poisson process on $(0,\infty)$ independent of an iid \seq\ 
$\bfQ_{i\cdot}=(\bfQ_{it})_{t\in\bbz}$, $i\in\bbz$, with the \ds\ $F_{\bfQ}$ of $\bfQ=(\bfQ_t)$. Then
$(\Gamma_i^{-1/\a},\bfQ_{i\cdot})$, $i=1,2,\ldots$, constitute the points of a Poisson process
on $(0,\infty)\times \ell^\a(\bbr^d)$ with intensity \ms\ given by 
$\mu(dy,d\bfq)=\a y^{-\a-1}\,dy\,F_\bfQ(d\bfq)$. We write $\bfW_j=\sum_{t\in\bbz} \bfQ_{jt}$ with generic element $\bfW$, $\|\bfQ_{j\cdot}\|_p^p=\sum_{t\in\bbz}|\bfQ_{jt}|^p$ and $\|\bfQ_{j\cdot}\|_\infty
=\max_{t\in\bbz}|\bfQ_{jt}|$. 

\bpr\label{lem:seriesrep}
Assume the conditions of Theorem~\ref{pr:prophybridch} and $\a<p$.
Then the following \pp\ \rep\ holds
\beam
\label{eq:series:rep}
 \big(\bfxi_\alpha,\zeta_{\alpha,p}\big) 
&\eqd& \Big(\sum_{j=1}^\infty \big(\Gamma_j^{-1/\alpha} \bfW_j-C_j^{(\alpha)} \E[\bfW]\big)\,,
\big\|\big(\Gamma_j^{-1/\alpha}\|\bfQ_{j\cdot}\|_p\big)_{j\ge 1}\big\|_p \Big)\,,
\eeam
where $C_i^{(\alpha)}=\tfrac{\alpha}{\red \alpha-1}\big(i^{(\alpha-1)/\alpha}-(i-1)^{(\alpha-1)/\alpha}\big) {\bf 1}_{(1,2)}(\alpha)$ and  $C_j^{(\alpha)} \E[\bfW]$ is interpreted as 
$\bf0$ for $\a\in (0,1)$.
\epr
The a.s. \con\  of the infinite series on the \rhs\
follows by the arguments in Section 1.4 of Samorodnitsky and
Taqqu~\cite{samorodnitsky:taqqu:1994} where one also uses the fact that 
$\E[|\bfW|^\a]<\infty$. The latter property also ensures that
$\max_{j\ge0}\Gamma_j^{-1/\a}\|\bfQ_{j\cdot}\|_\infty$ has a Fr\'echet \ds\ with index $\a$; see 
Example 2.5.4 in Mikosch and Wintenberger \cite{mikosch:wintenberger:2024}.
\begin{proof}
We start with an auxiliary result which will be the key to \eqref{eq:series:rep}.
Consider an iid \seq\ $(U_j)$ uniformly distributed on $(0,1)$ with generic element $U$  and the 
following quantities
\beao
\wt \bfS_n&:=& \sum_{j=1}^n \big(U_j^{-1/\a}\bfW_j-\E[U^{-1/\a}\,\bfW]\,\1_{(1,2)}(\a)\big)\,,\\
\wt \gamma_{n,p}^p&:=&\sum_{j=1}^n U_j^{-p/\a} \|\bfQ_{j\cdot}\|_p^p\,,\\
\wt \gamma_{n,\infty}&:=& \max_{j=1,\ldots,n} U_j^{-1/\a}\,\|\bfQ_{j\cdot}\|_\infty\,.
\eeao
We also assume that $(U_j)$ and $(\bfQ_{j\cdot})$ are independent. 
\ble\label{lem:aux1} Under the conditions of Proposition~\ref{lem:seriesrep} we have for $(\bfu,x,\la)\in\bbr^d\times \R_+^2$, and 
$\a\in (0,2)\backslash\{1\}$, $p>\a$,
\beao
\E \Big[
 \ex^{ i\,\bfu^\top (n^{-1/\a}\wt \bfS_n)-\lambda (n^{-p/\a}\wt \gamma_{n,p}^p)}\1(n^{-1/\a}\wt \gamma_{n,\infty}\le x)
\Big)
\Big]\to \Psi_\bfX(\bfu,x,\la)\,,\qquad \nto\,,
\eeao
where the right-hand is the limit \chdf-Laplace transform in Theorem~\ref{pr:prophybridch}.
\ele
\par\noindent
{{\em Proof of Lemma~\ref{lem:aux1}.}
Direct calculation yields
\beao
\lefteqn{\E \Big[
 \exp\Big( i\,\bfu^\top \dfrac{\wt S_n}{n^{1/\alpha}}-\lambda \,\dfrac{\wt \gamma_{n,p}^p}{n^{p/\alpha}}\Big)\,\1\big(n^{-1/\a}\wt \gamma_{n,\infty}\le x\big) 
\Big]} \\
& =& \ex^{-i\,\bfu^\top \E[U^{-1/\alpha} \bfW] \1_{(1,2)}(\alpha) n^{1-1/\alpha}}\\
&&\times \Big\{
\E \Big[\exp\Big(
i\,\bfu^\top\, \frac{U^{-1/\alpha}\,\bfW}{n^{1/\alpha}} -\lambda \,\frac{U^{-p/\alpha} \|\bfQ\|_p^p}{n^{p/\alpha}}
\Big)\,\1\big(n^{-1/\a} U^{-1/\a}\|\bfQ\|_\infty\le x\big)\Big]
\Big\}^n  \\
&= &\big(\ex^{-i\,\,n^{-1/\a}\bfu^\top \E[U^{-1/\alpha} \bfW]\,\1_{(1,2)}(\alpha) } \big )^n \\
&&\qquad  \times 
\Bigg\{\Bigg(
1+i\,n^{-1/ \a} \bfu^\top \E[U^{-1/\alpha} \bfW]\,\1_{(1,2)}(\alpha)  \\
&&\quad +\E
\int_1^\infty \Big(
\exp\big(i \bfu^\top\frac{y}{n^{1/\alpha}} \bfW -\lambda\, (\frac{y}{n^{1/\alpha}})^p \,
\|\bfQ\|^p_p\big)\,\1\big(y\,n^{-1/\a}\,\|\bfQ\|_\infty\le x\big)\\
&&\qquad -1- 
i \bfu^\top \frac{y}{n^{1/\alpha}}\,\bfW\,\1_{(1,2)}(\alpha) 
\big)  d(-y^{-\a}) \Big)\Bigg)
\Bigg\}^n.
\eeao
Changing variables $y /n^{1/\alpha}\mapsto z$, we obtain 
\beao
&& \big(\ex^{-i\,n^{-1/\a}\bfu^\top \E[U^{-1/\alpha} \bfW]\,\1_{(1,2)}(\alpha)}\big )^n \\
&& \times 
\Bigg\{
1+i n^{-1/\a}\,\bfu^\top \E[U^{-1/\alpha} \bfW]\,\1_{(1,2)}(\alpha) \\
&&\qquad +n^{-1} \E
\int_{n^{-1/\alpha}}^\infty \Big(
\exp\big(i\,z\bfu^\top \bfW-\lambda z^p \|\bfQ\|^p_p\big)\,\1(z\,\|\bfQ\|_\infty\le x)-1 \\  
&&\qquad - 
iz\,\bfu^\top\bfW\,\1_{(1,2)}(\alpha) 
\Big) d(-z^{-\a})\Bigg\}^n \\
&=& \Bigg\{ \Big( 1-i\,n^{-1/\a}\bfu^\top\E[U^{-1/\alpha} \bfW]\,\1_{(1,2)}(\alpha)+O(n^{-2/\alpha}) \Big ) \\
&& \times \Bigg( 1+ in^{-1/\a} \bfu^\top \E[U^{-1/\alpha} \bfW]\,\1_{(1,2)}(\alpha)\\
&&\qquad + n^{-1} \E \int_{\red n^{-1/\alpha}}^\infty \Big(
\exp\big(i z\bfu^\top \bfW-\lambda z^p \|\bfQ\|^p_p\big)\1(z\,\|\bfQ\|_\infty\le x)-1\\
&&\qquad - 
i\,z\,\bfu^\top \bfW \,\1_{(1,2)}(\alpha) 
\Big) d(-z^{-a})  \Bigg) \Bigg \}^n \\
&=& \Bigg\{
1+n^{-1} \int_{\red n^{-1/\alpha}}^\infty \E\big[
\exp\big(iz\,\bfu^\top \bfW-\lambda z^p \|\bfQ\|^p_p\big)\,\1(z\,\|\bfQ\|_\infty\le x)-1 \\
&&\qquad - 
i\,z\,\bfu^\top \bfW\, \1_{(1,2)}(\alpha) 
\big) 
\big]d (-z^{-\alpha}) +o(n^{-1}) 
\Bigg\}^n\,,
\eeao
where the change of integral and expectation in the penultimate step
 is justified by Fubini. 
Passing $n\to \infty$, we obtain the desired limit $\Psi(\bfu,x,\la)$
for every $(\bfu,x,\la)\in\bbr^d\times \R_+^2$. \qed
\par
We conclude from Lemma~\ref{lem:aux1} that
\beao
n^{-1/\a}\big(\wt \bfS_n\,,\wt \gamma_{n,p}\,,\wt \gamma_{n,\infty}\big)
\std (\bfxi_\a,\zeta_{\a,p},\zeta_{\a,\infty})\,,\qquad \nto\,.
\eeao
\par
We recall that, given $\Gamma_{n+1}$, the vector 
$(\Gamma_i/\Gamma_{n+1})_{i=1,2,\ldots,n}$ has the same \ds\  
as the order statistics of the sample $U_1,\ldots,U_n$.
We restrict ourselves to the case $\a\in (1,2)$, the case $\a\in (0,1)$ 
being similar. We have

\beao
\lefteqn{\Big(
\dfrac{\wt \bfS_n}{n^{1/\alpha}}, \dfrac{\wt \gamma_{n,p}^p}{n^{p/\alpha}}\,,\dfrac{\wt \gamma_{n,\infty}}
{n^{1/\a}}
\Big) }\\ 
&\eqd&\Big(
\Big(\dfrac{\Gamma_{n+1}}{n} \Big)^{1/\alpha} \sum_{i=1}^n \Gamma_i^{-1/\alpha} \bfW_i
-\dfrac{\alpha}{\red \alpha-1} n^{1-1/\alpha}\E [\bfW],\, \Big(\dfrac{\Gamma_{n+1}}{n}\Big)^{p/\alpha} \sum_{i=1}^n 
\Gamma_i^{-p/\alpha} \|\bfQ_{i\cdot}\|_p^p\,,\\
&& \Big(\dfrac{\Gamma_{n+1}}{n}\Big)^{1/\a}\, \max_{i=1,\ldots,n} \Gamma_i^{-1/\a}\,\|\bfQ_{i\cdot}\|_\infty  \Big)\,.
\eeao

Since $\Gamma_{n+1}/n \to 1$ and $\Gamma_i>0$ a.s. the first component on the 
\rhs\ converges a.s. to the desired limit:

\beao
\lefteqn{ \Big( \dfrac{\Gamma_{n+1}}{n}
\Big)^{1/\alpha} 
\sum_{i=1}^n \Gamma_i^{-1/\alpha} \bfW_i - 
\dfrac{\alpha}{\alpha-1} n^{1-1/\alpha}\E [\bfW]} \\
& =& \Big(
\dfrac{\Gamma_{n+1}}{n}
\Big)^{1/\alpha} \sum_{i=1}^n \Big(
\Gamma_i^{-1/\alpha} \bfW_i - \E [\bfW] \dfrac{\alpha}{\alpha-1}\big(
i^{(\alpha-1)/\alpha}-(i-1)^{(\alpha-1)/\alpha}
\big)
\Big) \\
& & + \dfrac{\alpha}{\alpha-1} n^{(\alpha-1)/\alpha} \big(
(n^{-1}\Gamma_{n+1})^{1/\alpha} -1
\big) \E [\bfW] \\
& \stas &\sum_{i=1}^\infty \big(
\Gamma_i^{-1/\alpha} \bfW_i - C_i^{(\alpha)} \E[\bfW]
\big),\qquad \red \nto\,. 
\eeao
{\red Indeed, by Taylor expansion with \pro y~1\,,
\beao\lefteqn{
n^{(\alpha-1)/\alpha} \big(
(n^{-1}\Gamma_{n+1} \pm 1)^{1/\alpha} -1
\big)}\\ 
&=& n^{(\alpha-1)/\alpha} \big(\Gamma_{n+1}/n-1\big) O(1)\\
&=&n^{-1/\a}\big(\Gamma_{n+1}-(n+1)\big)\, O(1)=o(1)\,,\qquad \nto\,.
\eeao
In the last step we used the law of the iterated logarithm and the fact that $\a<2$.
}
The condition $p/\alpha>1$ implies
\beao
 \Big(\dfrac{\Gamma_{n+1}}{n}
\Big)^{p/\alpha} \sum_{i=1}^n \Gamma_i^{-p/\alpha} \|\bfQ_{i\cdot}\|_p^p \stas
\sum_{i=1}^\infty \Gamma_i^{-p/\alpha} \|\bfQ_{i\cdot}\|_p^p\,,\qquad \nto\,, 
\eeao
and we also have 
\beao
\Big(\dfrac{\Gamma_{n+1}}{n}\Big)^{1/\a}\,\max_{i=1,\ldots,n} 
\Gamma_i^{-1/\a}\,\|\bfQ_{i\cdot}\|_\infty \stas \max_{i\ge 1} 
\Gamma_i^{-1/\a}\,\|\bfQ_{i\cdot}\|_\infty \qquad \nto\,.
\eeao
Thus we have proved the proposition.} \end{proof}
In what follows, we will conduct a change of \ms\ on $\bfQ$ for $p>\a$:
\beam\label{eq:tildeQ}
\P ( \bfQ^{(p)} \in \cdot) := \E\Big[\dfrac{\|\bfQ\|_p^\alpha}{\E[\|\bfQ\|_p^\alpha]}
 \1\Big(\dfrac{\bfQ}{\|\bfQ\|_p}\in \cdot\Big)\Big]\,;
\eeam
the case $p=\infty$ corresponds to $\|\bfQ\|_\infty=\max_{t\in\bbz} |\bfQ_t|$. We also observe that 
the case $\a=p$ is formally included: $\bfQ^{(\a)}\eqd \bfQ$ since $\|\bfQ\|_\alpha=1$ a.s.
\par
Suppressing the dependence on $\a$, we write 
\beam\label{eq:suppress}
 \bfW^{(p)}:=
\sum_{t\in\bbz} \bfQ_t^{(p)}\,,\qquad p>\a\,.
\eeam
\par 
By a change of measure we get an alternative expression for the \ds\ of 
$(\bfxi_\alpha,\zeta_{\alpha,p})$ in terms of their \chf-Laplace transform.
\ble\label{prop:Wtilde} Under the conditions of Proposition~\ref{lem:seriesrep}
we have for $(\bfu,\la)\in\bbr^d\times \R_+$, and $\a\in (0,2)\backslash\{1\}$, $p>\a$,
\beam\label{j:ch:LP}
\lefteqn{\E\big[\exp\big(i\bfu^\top\bfxi_\alpha-\lambda \,\zeta_{\alpha,p}^p\big)\big]}
\nonumber\\
&=& \exp\Big(\E[\|\bfQ\|_p^\alpha]\, 
\int_0^\infty \E\Big[\exp\Big(i\,y\,\bfu^\top  \bfW^{(p)} -\lambda y^p\Big)-1-i\,y\,
\bfu^\top  \bfW^{(p)} \,\1_{(1,2)}(\alpha)\Big] d(-y^{-\alpha})\Big)\,.\nonumber
\eeam
\ele
\begin{proof}
Writing $\Psi_\bfX(\bfu,\la):=\Psi_\bfX(\bfu,\infty,\la)$ and 
changing variables $y\,\|\bfQ\|_p\mapsto z$ \textcolor{red}{after Fubini}, 
we obtain
\beao
&& \log \Psi_\bfX(\bfu,\la) \\
&=&
\int_0^\infty \E \Big[\exp\big(iy\,\bfu^\top \bfW-y^p\,\lambda \|\bfQ\|_p^p\big)-1-i
y\,{\red \bfu^\top}\,\bfW \1_{(1,2)}(\alpha)\Big]\, d(-y^{-\alpha})\\
&=& \textcolor{red}{\E \Big[} \int_0^\infty  \|\bfQ\|_p^\alpha\,  
\Big( \exp\big(iz\,\bfu^\top  \bfW/\|\bfQ\|_p-\lambda\, z^p\big)-1-i\,z\,\bfu^\top  
\bfW/\|\bfQ\|_p \,\1_{(1,2)}(\alpha)\Big) d(-z^{-\alpha})\Big] \\
& =& \E[\|\bfQ\|_p^\alpha] \int_0^\infty \E\Big[\exp\big( i\,z\,\bfu^\top  \bfW^{(p)}-
\lambda\, z^p\big)-1-iz\,\bfu^\top  \bfW^{(p)} \,\1_{(1,2)}(\alpha)\textcolor{red}{\Big]} \,d(-z^{-\alpha})\,.
\eeao
In view of Proposition~\ref{lem:seriesrep} this is the logarithm of the joint 
\chf-Laplace transform of  $(\bfxi_\a,\zeta_{\a,p}^p)$. 
\end{proof}
\section{Conditions for finite moments of $\bfR_{\a,p}$}\label{sec:mom}\setcounter{equation}{0}
Denote an iid \seq\ with generic element $ \bfW^{(p)}$ by $( \bfW_i^{(p)})$
and assume that it is independent of the Poisson points $(\Gamma_i)$. 
Using  Lemma~\ref{j:ch:LP}, we can derive 
series \rep s which are analogous to Proposition~\ref{lem:seriesrep}:
\begin{align}
\label{measure:change:series:expression}
 \big(\bfxi_\alpha,\zeta_{\alpha,p}) \eqd 
\big(\E[\|\bfQ\|_p^\alpha]\big)^{1/\alpha}\Big( \sum_{i=1}^\infty \Big( \Gamma_i^{-1/\alpha} 
\bfW_i^{(p)}-C_i^{(\alpha)}\, \E[ \bfW^{(p)}]\Big),\,
\big\|\big(\Gamma_i^{-1/\alpha}\big)_{i\ge 1}\big\|_p \Big)\,.
\end{align}
An immediate con\seq\ of \eqref{measure:change:series:expression}
is that for $p>\a$,
\beao
\bfR_{\a,p}:=\dfrac{\bfxi_\a}{\zeta_{\a,p}}&\eqd &
\dfrac{\sum_{i=1}^\infty \big( \Gamma_i^{-1/\alpha}  \bfW_i^{(p)}-
C_i^{(\alpha)}\, \E[ \bfW^{(p)}]\big)}{\big\|\big(\Gamma_i^{-1/\alpha}\big)_{i\ge 1}\big\|_p }\\
&=& \underbrace{\dfrac{\sum_{i=1}^\infty \big( \Gamma_i^{-1/\alpha}
\bfW_i^{(p)}-C_i^{(\alpha)}\, \E[ \bfW^{(p)}]\big)}{\big\|\big(\Gamma_i^{-1/\alpha}|\bfW^{(p)}_i|\big)_{i\ge 1}\big\|_p }}_{=:\rho_1} 
\underbrace{ 
\dfrac{\big\|\big(\Gamma_i^{-1/\alpha}|\bfW^{(p)}_i|\big)_{i\ge 1}\big\|_p}
{\big\|\big(\Gamma_i^{-1/\alpha}\big)_{i\ge 1}\big\|_p }}_{=:\rho_2}\,.
\eeao
We observe that the condition $\E[|\bfW^{(p)}|^\a]<\infty$, $p>\a$, holds under our assumptions since
$\bfW=\bfW^{(\a)}$ and $ \E[|\bfW|^\a]<\infty$ in view of Remark \ref{rem:moments}. Hence 
Theorem 1.4.5 in \cite{samorodnitsky:taqqu:1994} applies to the series \rep\ of $\bfxi_\a$ in \eqref{measure:change:series:expression}, ensuring 
the a.s. \con\ of this series and of $\rho_1$. 
The ratio $\rho_1$ corresponds to the iid situation described in Section~\ref{subsec:iid} with non-normalized spectral component $\bfW^{(p)}$. In this case, 
$\rho_1$ has all moments finite. The factor $\rho_2$ has $q$th moment for $q\ge p$ if $\E[|\bfW^{(p)}|^q]<\infty$. Indeed, writing $p_i:=\Gamma_i^{-p/\a}/\|(\Gamma_j^{-1/\a})_{j\ge 1}\|_p^p$ and observing that  $\sum_{i\ge 1}p_i=1$ a.s., we have by Lyapunov's inequality
\beao
\E[\rho_2^q]&=&\E\Big[ \Big(\sum_{i\ge 1}^\infty p_i\,|\bfW_i^{(p)}|^p \Big)^{q/p}\Big]
\le \E\Big[\sum_{i\ge 1}^\infty p_i\,|\bfW_i^{(p)}|^q \Big]=\E\big[|\bfW^{(p)}|^q\big]\,.
\eeao
Applying H\"older's inequality, we obtain for $q\ge p>\a$ and $q_1,q_2>0$ such that $q_1^{-1}+q_2^{-1}=1$,
\beao
\E\big[|\bfR_{\a,p}|^q\big]\le \big(\E\big[|\rho_1|^{q\,q_1}\big]\big)^{1/q_1}\,\big(\E\big[\rho_2^{q\,q_2}\big]\big)^{1/q_2}\,.
\eeao
Since we can choose $q_2$ arbitrarily close to 1 we conclude that $\bfR_{\a,p}$ has $q$th moment for $q\ge p>\a$ if $\E[|\bfW^{(p)}|^{q+\vep}]<\infty$ for some $\vep>0$.
For integer-valued
$q$ we can make even more precise statement about the relation of moments between $\bfR_{\a,p}$ and $\bfW^{(p)}$; see Theorem~\ref{thm:momentsratio}
below.
\par
Another interesting con\seq\ of \eqref{measure:change:series:expression}
is that for $p>\a$, the limit ratio $|\bfR_{\a,p}|$ is stochastically dominated by the quantity
\begin{align}
\label{eq:dominfty}
\dfrac{\big|\sum_{i=1}^\infty \big( \Gamma_i^{-1/\alpha}  \bfW_i^{(p)}-
C_i^{(\alpha)}\, \E[ \bfW^{(p)}]\big)\big|}{\big\|\big(\Gamma_i^{-1/\alpha}\big)_{i\ge 1}\big\|_\infty }=\dfrac{\big|\sum_{i=1}^\infty \big( \Gamma_i^{-1/\alpha}  \bfW_i^{(p)}-
C_i^{(\alpha)}\, \E[ \bfW^{(p)}]\big)\big|}{\Gamma_1^{-1/\alpha}}.
\end{align}
Thus, results about the existence of moment of $\bfR_{\a,p}$ can be reduced to the corresponding results for \eqref{eq:dominfty}.
This is the content of the following result.
\bth\label{thm:momentsratio}
 Consider a   \regvary\ stationary \seq\ $(\bfX_t)$ satisfying the conditions of 
Theorem~\ref{pr:prophybridch} and $p\in (\a,\infty]$. 
\begin{enumerate}
\item[\rm 1.]
Assume $\a\in  (0,1)$. Then  for every integer $m\ge 1$, 
$\E[|\bfW^{(p)}|^m]<\infty$ implies that $\E[|\bfR_{\a,p}|^m]<\infty$.
\item[\rm 2.] Assume
$\a\in  (1,2)$. Then  for every even integer $m\ge 2$, $\E[|\bfW^{(p)}|^m]<\infty$ 
implies that $\E[|\bfR_{\a,p}|^m]<\infty$. 
\item[\rm 3.]
If $p=\infty$ the sufficient moment conditions in 1. and 2. are also necessary for $\E[|\bfR_{\a,p}|^m]<\infty$.
\end{enumerate}
\ethe
\begin{proof}
Without loss of generality we may restrict ourselves to the case $d=1$. 
\subsection*{The case $\a\in (0,1)$} In absence of centering in \eqref{eq:dominfty}, we  immediately get
\beao
\dfrac{\Big|\sum_{i=1}^\infty   \Gamma_i^{-1/\alpha}  W_i^{(p)}\Big|}{\big\|\big(\Gamma_i^{-1/\alpha}\big)_{i\ge 1}\big\|_p }
\le\dfrac{\sum_{i=1}^\infty   \Gamma_i^{-1/\alpha}  |W_i^{(p)}|}{\big\|\big(\Gamma_i^{-1/\alpha}\big)_{i\ge 1}\big\|_\infty }=:R_\a^{(p)}\,,
\eeao
where $W^{(p)}$ is the univariate equivalent of $\bfW^{(p)}$ in \eqref{eq:suppress}.
\par
The proof of 1. and 3. follows from this inequality and the following lemma.

\ble\label{lem:x} Consider a real-valued \regvary\ stationary \seq\ $(X_t)$ with index $\a\in (0,1)$ satisfying the conditions of 
Theorem~\ref{pr:prophybridch}. Then for every integer $m\ge 1$, $\E[|W^{(p)}|^m]<\infty$ \fif\ $\E[|R_\a^{(p)}|^m]<\infty$.
\ele

\begin{proof}[Proof of Lemma~\ref{lem:x}] 
Matsui et al. \cite{matsui:mikosch:wintenberger:2024}, Corollary~4.1,  derived the Laplace transform of $R_\a^{(p)}$:
\beao
  \Psi_{R_\a^{(p)}}(\lambda) = \dfrac{\E\big[\exp\big(-\lambda \,|W^{(p)}|\big)\big]}{
\int_0^\infty \E \big[1-\exp\big(-y\, \lambda \,|W^{(p)}|\big) \,\1(y\le 1)\big]\,d(-y^{-\alpha})}=:
\dfrac{h(\la)}{g(\la)}\,,\qquad \la>0\,.
\eeao
\par
First we assume that  $\E[|W^{(p)}|^m]<\infty$. 
By monotone \con , $\lim_{\lambda\downarrow 0}g(\lambda)=1$. Calculation yields the $m$th derivative of $g$:
\beao 
 g^{(m)}(\lambda)=(-1)^{m-1} \int_0^1 \E\big[|W^{(p)}|^m \exp\big(-\lambda y\,|W^{(p)}|\big)\big] \,\alpha y^{m-\alpha-1}\, dy\,,
\eeao
and monotone convergence implies
\beao
 \lim_{\la\downarrow   0} |g^{(m)}(\lambda) |= 
\dfrac{\alpha}{m-\alpha}\, \E[|W^{(p)}|^m]<\infty.
\eeao
Obviously, $\lim_{\la \downarrow 0}\E\big[|W^{(p)}|^m \,\exp\big(-\lambda \,|W^{(p)}|\big)\big]=\E[|W^{(p)}|^m]<\infty$.
\par
We observe that 
\beam
\label{derivative:phi:m}
 \Psi_{R_\a^{(p)}}^{(m)}(\la)= \sum_{k=0}^m {m \choose k}\, 
h^{(k)}(\la)
(g^{-1}(\lambda))^{(m-k)}\,,
\eeam
where $g^{-1}=1/g$.
Since $(g^{-1})^{(n)}$\,, $n\le m$, are linear combinations of products composed by $g^{(\ell)}/g^n,\,\ell,n\le m$, 
\beao
 \E[|R_\a^{(p)}|^m]=(-1)^m \,\Psi_{R_\a^{(p)}}^{(m)}(0) =(-1)^m \,\lim_{\lambda \downarrow 0}\Psi_{R_\a^{(p)}}^{(m)}(\lambda)<\infty. 
\eeao
\par
Now assume that
$\E[|R_\a^{(p)}|^m]<\infty$ for some integer $m\ge 1$. This relation implies that 
$\E[|R_\a^{(p)}|^k]=(-1)^k \,\Psi_{R_\a^{(p)}}^{(k)}(0)<\infty$ 
for $k\le m$. Observe that 
\beao
 -\Psi_{R_\a^{(p)}}^{(1)}(\lambda)&=& \dfrac{\E[|W^{(p)}|\,\exp\big(-\lambda\,|W^{(p)}|\big)\big]}{g(\lambda)}\\&& + 
\dfrac{h(\la)}{g^2(\lambda)}
\,\int_0^1 \E \big[|W^{(p)}|\,\exp\big(-\lambda\, |W^{(p)}|\big)\big]\,\a\, y^{-\alpha} \,dy\,.
\eeao
By monotone convergence,
\beao
 \lim_{\la\downarrow 0} \int_0^1 \E\big[|W^{(p)}| \,\exp\big(-\lambda\, |W^{(p)}|\big)\big]\, \alpha\, y^{-\alpha}\,dy =
\dfrac{\alpha}{1-\alpha}\, \E[|W^{(p)}|]\,. 
\eeao
Thus, letting $\lambda \downarrow 0$, we obtain 
\beao
 -\Psi_{R_\a^{(p)}}^{(1)}(0) = \dfrac{\E[|W^{(p)}|]}{1-\alpha}<\infty\,. 
\eeao
\par
Now we proceed by induction on $k\le m$. Assume that
$\E[|R_\a^{(p)}|^k]=(-1)^k \,\Psi_{R_\a^{(p)}}^{(k)}(0)<\infty$ implies $\E[|W^{(p)}|^k]<\infty$. We will show that
$1\le k<m$ may be replaced by $k+1$. In view of \eqref{derivative:phi:m} we have 
\beao
 \Psi_{R_\a^{(p)}}^{(k+1)}(\lambda) &=& \sum_{l=0}^{k+1} {k+1 \choose l} 
h^{(l)}(\la)\,
(g^{-1}(\lambda))^{(k+1-l)} \\
&=&h^{(k+1)}
g^{-1}(\la) + (g^{-1}(\lambda))^{(k+1)}\,h(\la)
 +G_k(\lambda)\,,
\eeao
where $G_k(\lambda)$ is a linear combination of $\big(\E\big[\exp\big(-\lambda\, |W^{(p)}|\big)\big]\big)^{(l)}\,(g^{-1}(\lambda))^{(n)},\,l,n\le k$.
By the induction hypothesis,  $\lim_{\lambda\downarrow 0}G_k(\lambda)=G_k(0)$ exists and is finite. 
Thus we consider the limits of the first two terms in the last expression. 
Clearly, 
\beao
h^{(k+1)}(\la)
= (-1)^{k+1} 
\E\big[|W^{(p)}|^{k+1}\exp\big(-\lambda\, |W^{(p)}|\big)\big]\,,
\eeao
and the second term is given by
\beao
&& -\sum_{l=0}^k \binom{k}{l} g^{(l+1)}(\lambda)\,(g^{-2}(\lambda))^{(k-l)}\,h(\la)
= \big(
-g^{(k+1)}(\lambda) g^{-2}(\lambda) + H_k(\lambda)\big)\,h(\la)\,, 
\eeao
where $H_k(\lambda)$ is a linear combination of products composed by $g^{(l)}(\lambda),\,l\le k$, and powers of $g^{-1}(\lambda)$.
Hence $\lim_{\lambda \downarrow 0} H_k(\lambda)$ exists. We observe that 
\[
 -g^{(k+1)}(\lambda)g^{-2}(\lambda) = (-1)^{k+1} \int_0^{1}
\E\big[|W^{(p)}|^{k+1}\exp\big(-\lambda y\, |W^{(p)}|\big)\big] \alpha y^{k-\alpha} dy \, g^{-2}(\lambda).  
\]
We conclude that $\lim_{\lambda \downarrow 0} \Psi_{R_\a^{(p)}}^{(k+1)}(\lambda)$ is a constant plus the limit of 
\beao
& & h^{(k+1)}(\lambda)\,
g^{-1}(\lambda) -g^{(k+1)}(\lambda) g^{-2}(\lambda)\,h(\la) \\
& =& (-1)^{k+1} g^{-1}(\lambda) \int_0^1
\E\big[|W^{(p)}|^{k+1} \big(
\exp\big(-\lambda \, |W^{(p)}|\big)+ \exp\big(-\lambda y\, |W^{(p)}|\big) \alpha y^{k-\alpha} 
\Psi_{R_\a^{(p)}}(\lambda)
\big) \big ] dy \,.
\eeao
The integrand increases as $\lambda \downarrow 0$ since $\Psi_{R_\a^{(p)}}(\lambda)$ is a Laplace transform. 
Finally, we obtain for some constant $c$,
\beao
 \Psi_{R_\a^{(p)}}^{(k+1)}(0)= (-1)^{k+1}\, \E[|W^{(p)}|^{k+1}]\dfrac{k+1}{k+1-\alpha} +c\,.
\eeao
Thus we have proved $\E[|W^{(p)}|^{k+1}]<\infty$ for $1\le k<m$.
\end{proof}
\subsection*{The case $\a\in (1,2)$} Define
\beao
R_\a^{(p)}:=\dfrac{ \sum_{i=1}^\infty \big( \Gamma_i^{-1/\alpha} W_i^{(p)}-
C_i^{(\alpha)}\, \E[W^{(p)}]\big) }{\big\|\big(\Gamma_i^{-1/\alpha}\big)_{i\ge 1}\big\|_\infty }
\eeao
and observe that $|R_{\a,p}|\le |R_\a^{(p)}|$ with identity for $p=\infty$. 
 In the following result we provide conditions for the existence of even $m$th moment  $\E[|R_\alpha^{(p)}|^m]$ in terms of the moments of $W^{(p)}$.
These imply the existence of the $m$th moment of $R_{\a,p}$ for $p\in (\a,\infty]$ and also give necessary conditions for 
the existence of $\E[R_{\a,\infty}^m]$. In turn, 2. and 3. are proved for $\a\in (1,2)$.
\ble\label{lem:2} Consider a real-valued \regvary\ stationary \seq\ $(X_t)$ with index $\a\in (1,2)$ and $\E[X]=0$, satisfying the conditions of 
Theorem~\ref{pr:prophybridch}, and choose an even integer $m\ge 2$. Then 
$\E[(R_\alpha^{(p)})^m]<\infty$ holds if and only if $\E[(W^{(p)})^{m}]<\infty$. 
\ele
\begin{proof} We appeal to the \chf\ of $R_\a^{(p)}$ in Matsui et al. \cite{matsui:mikosch:wintenberger:2024}, Corollary~4.1:
\beam\label{chfalpha2}
 \varphi_{R_\alpha^{(p)}}(u)&=& \dfrac{\E[\exp\big(i\,u \,W^{(p)}\big)]}{
\int_0^\infty \E \big[1+i\,y\,u \,W^{(p)}-\exp\big(i\,y\,u W^{(p)}\big)\,\1_{[0,1]}(y)\big] \,d(-y^{-\alpha})}=:\dfrac{f(u)}{g(u)}\,,\qquad u\in\bbr\,. \nonumber\\
\eeam
The existence of  $\E[|W^{(p)}|]$ follows from Remark 4.2 in \cite{matsui:mikosch:wintenberger:2024}.
By virtue of \eqref{chfalpha2} we obtain
\beao
 \varphi_{R_\alpha^{(p)}}'(u)= \varphi_{R_\a^{(p)}}(u)\,
\Big(\dfrac{f'(u)}{f(u)}-\dfrac{g'(u)}{g(u)}\Big)\,,\qquad u\in \bbr\,. 
\eeao
Moreover, under the condition $\E[(W^{(p)})^m]<\infty$ for an even $m$, the $m$th derivatives of $f,g$ at $u \neq 0$ 
are well-defined since  
\begin{align*}
 |g^{(m)}(u)| &= \Big|
i^{m-2} \int_0^1 \E[ (W^{(p)} )^m] \ex^{iy u W^{(p)}} y^m \,d(-y^{-\alpha})
\Big| \\
&\le \int_0^1 \E[  (W^{(p)} )^m] y^m d(-y^{-\alpha}) <\infty,\quad u \neq 0\,,
\end{align*}
and $|f^{(m)}(u)|<\infty$ holds in a similar way. 
Then, differentiation of $\varphi'_{R_\alpha^{(p)}}$ further yields for $m\ge 1$,
\beam
\label{rep:Ralpham}
\varphi_{R_\alpha^{(p)}}^{(m)}(u) = \sum_{k=0}^{m-1} \binom{m-1}{k} 
\varphi_{R_\a^{(p)}}^{(k)}(u)\, \Big(
\dfrac{f'(u)}{f(u)}-\dfrac{g'(u)}{g(u)}\Big)^{(m-1-k)}\,,\qquad u\in\bbr\,,
\eeam
where $(\cdot)^{(k)}$ denote $k$th differentiation of a function inside the parentheses.
\par
Assuming $\E[(W^{(p)})^m]<\infty$ for some even $m$,
we will show $\E[(R_\alpha^{(p)})^m]<\infty$ by proving the existence of $\lim_{u\to 0}\varphi_{R_\alpha^{(p)}}^{(m)}(u)$. 
First we will show that the limits
\beam\label{lim:gg:ff}
 \lim_{u\to 0}\Big(\dfrac{f'(u)}{f(u)}\Big)^{(m-1)}\qquad \text{and}\qquad \lim_{u\to 0}\Big(\dfrac{g'(u)}{g(u)}\Big)^{(m-1)} 
\eeam
exist. We restrict ourselves to prove \eqref{lim:gg:ff} for $g$; the $f$-case is similar. Direct calculation yields
 $g(0)=1$, $g'(0)=i(\alpha/(\alpha-1))\,\E[W^{(p)}]$, {\red and $g''(0)=(\a/(\a-2))\E[(W^{(p)})^2]$}. 
Hence \eqref{lim:gg:ff} holds for $g$ and $m=2$. Now assume that $m>2$ is even. 
We proceed by induction: assume that $\lim_{u\to 0}(g'(u)/g(u))^{(j)}$ exists for $j\le m-2$ and prove its existence for $j=m-1$. 
In view of the expression 
\beao
 \Big(\dfrac{g'(u)}{g(u)}\Big)^{(m-1)} &= &\sum_{k=0}^{m-1} \binom{m-1}{k} g^{(1+k)}(u) \,\big(g^{-1}(u)\big)^{(m-1-k)}, 
\eeao
it suffices to observe that 
\beao
 \big|
\lim_{u\to 0} g^{(m)}(u)
\big| = \big| -i^m \dfrac{\alpha}{m-\alpha} \E[(W^{(p)})^m] \big|= \dfrac{\alpha}{m-\alpha} \E[(W^{(p)})^m]<\infty. 
\eeao
Next we proceed to the limit of \eqref{rep:Ralpham}. Again we use the induction hypothesis. Suppose that the following limit exists:
\beao
 \varphi_{R_\a^{(p)}}^{(j)}(0) := \lim_{u\to 0} \varphi_{R_\a^{(p)}}^{(j)}(u)\,,\qquad j\le m-2\,,
\eeao
and we will show that $\varphi_{R_\a^{(p)}}^{(m)}(0)$ exists. For $j=1$ this holds by virtue of Remark 4.2 in \cite{matsui:mikosch:wintenberger:2024}.
In view of  \eqref{rep:Ralpham} and \eqref{lim:gg:ff} it suffices to check the existence of 
$\lim_{u\to 0} \varphi_{R_\a^{(p)}}^{(m-1)}(u)$. Replacing $m$ by $m-1$ in \eqref{rep:Ralpham} and  using  
the existence of the limits \eqref{lim:gg:ff} and $\varphi_{R_\a^{(p)}}^{(j)}(0)$\,, $j\le m-2$,  we can prove the existence of 
$\varphi_{R_\alpha^{(p)}}^{(m)}(0)$ which implies $\E[(R_\alpha^{(p)})^m]<\infty$. 
\par
Now assume that  $\E[|W^{(p)}|^{m-1}]<\infty$,  $ \E[(W^{(p)})^m]=\infty$ for some even $m$. We will prove that $ \E[(R_\a^{(p)})^m]=\infty$
by showing that the derivative \eqref{rep:Ralpham} is not finite. We observe that $\varphi_{R_\a^{(p)}}^{(k)}(u)$\,, $k\le m-1$, exist and are  finite, and so do
\beao
 \lim_{u\to 0} \Big(
\dfrac{f'(u)}{f(u)}-\dfrac{g'(u)}{g(u)}
\Big)^{(k)},\qquad k\le m-2\,.
\eeao
Thus it suffices to check the term in \eqref{rep:Ralpham} corresponding to $k=0$. The limit can be written in the form
\beao
 \lim_{u\to 0} \Big(
\dfrac{f'(u)}{f(u)}-\dfrac{g'(u)}{g(u)}\Big)^{(m-1)}= \lim_{u\to 0}\Big(
\dfrac{f^{(m)}(u)}{f(u)}-\dfrac{g^{(m)}(u)}{g(u)}
\Big) +c\,, 
\eeao
for some constant $c$.
However, the limit does not exist. If it existed we would have  
\beao
\lim_{u\to 0} \Big|
\dfrac{f^{(m)}(u)}{f(u)}-\dfrac{g^{(m)}(u)}{g(u)}
\Big| = \Big| 
\E [(W^{(p)})^m]\,i^m\, \Big(1+ \dfrac\alpha { m-\alpha}\Big)
\Big|<\infty, 
\eeao
in contradiction to our assumption. 
In a similar manner we can prove $ \E[(R_\a^{(p)})^m]=\infty$ under the assumption that 
$\E[(W^{(p)})^{m-2}]<\infty, \E[|W^{(p)}|^{m-1}]=\infty$ but 
$ \Big|
\dfrac{f^{(m-1)}(0)}{f(0)}-\dfrac{g^{(m-1)}(0)}{g(0)}
\Big|<\infty$ for some even $m$.
Finally, we assume that $\E[(W^{(p)})^{m-2}]<\infty$ and 
$ \Big|
\dfrac{f^{(m-1)}(0)}{f(0)}-\dfrac{g^{(m-1)}(0)}{g(0)}
\Big|=\infty$ for some even $m$.
In view of \eqref{rep:Ralpham} we have 
\[
 \varphi_{R_\alpha^{(p)}}^{(m-1)}(u) = \varphi_{R_\alpha^{(p)}}(u) 
\Big( \dfrac{f^{(m-1)}(u)}{f(u)}-\dfrac{g^{(m-1)}(u)}{g(u)} \Big)+h(u)\,,\qquad u \neq 0\,,
\]
for some function $h(u)$ such that $\lim_{u\to\infty}h(u)$ exists and is finite.
Therefore $\lim_{u\to0}\varphi_{R_\alpha^{(p)}}^{(m-1)}(u)$ cannot exist neither $\lim_{u\to0}\varphi_{R_\alpha^{(p)}}^{(m)}(u)$. The latter fact is equivalent to $\E[(R_\a^{(p)})^m]=\infty$ because $m$ is even.
\end{proof}
{\red This finishes the poof of Theorem~\ref{thm:momentsratio}.}
\end{proof}
Applications of Theorem~\ref{thm:momentsratio} require to understand the moment properties of $W^{(p)}$. 
In what follows, we will provide sufficient conditions for the existence of $\E[(W^{(p)})^m]$ in terms of the 
forward spectral tail process.
\par
By \eqref{eq:tildeQ} we observe that
\beam\label{eq:finite}
\E\big[(W^{(p)})^m]&\le & \E \big[\|Q\|_p^{\a-m}  |\textcolor{red}{W}|^m \big] /\E\big[\|Q\|_p^\a\big]
\le\E \big[\|Q\|_p^{\a- m} \, \|Q\|_1^m \big] /\E\big[\|Q\|_p^\a\big]\,.\nonumber\\
\eeam
{\red It suffices to check that the numerator on the \rhs\ is finite. For this reason,
recall the time-change formula (see Basrak and Segers \cite{basrak:segers:2009}; cf.  Section 3.3 of \cite{mikosch:wintenberger:2024}): for every $t\in \Z$,
\beam\label{eq:tcpds}
\P\big((\Theta_s)_{s\in \Z} \in \cdot , \Theta_{t}\ne 0\big)&=&\E\Big[|\Theta_{-t}|^\alpha \Big(\dfrac{(\Theta_{s-t})_{s\in \Z}}{|\Theta_{-t}|}\in \cdot\Big)\Big]\,.
\eeam
We write $\Theta_{s:\infty}:=(\Theta_t)_{t\ge s}$, $s\in\bbz$, and apply  \eqref{eq:tcpds} 
to the right-hand numerator in \eqref{eq:finite}. We obtain
an expression only in terms of the forward process $(\Theta_t)_{t\ge 0}$:}
\beao
\E\big[\|Q\|_1^m\,\|Q\|_p^{\alpha-m}\big]
&=&\E\big[\|\Theta\|_1^m\,\|\Theta\|_p^{\alpha-m}/\|\Theta\|_\alpha^\alpha\big]\\
&=&\sum_{t\in\Z}\E\big[(\|\Theta_{t:\infty}\|_1^m-\|\Theta_{t+1:\infty}\|_1^m) / (\|\Theta\|_p^{m-\alpha}\|\Theta\|_\alpha^\alpha)] \\
&=&\sum_{t\in\Z}\E\big[|\Theta_{-t}|^\alpha (\|\Theta_{0:\infty}\|_1^m - \|\Theta_{1:\infty}\|_1^m)/(\textcolor{red}{
\|\Theta\|_p^{m-\alpha}}\|\Theta\|_\alpha^\alpha)] \\
&=&\E\big[\|\Theta\|_\alpha^\alpha(\|\Theta_{0:\infty}\|_1^m-\|\Theta_{1:\infty}\|_1^m)/(\|\Theta\|_p^{m-\a}\,\|\Theta\|_\alpha^\alpha)\big]\\
&=&\E\big[ (\|\Theta_{0:\infty}\|_1^m -\|\Theta_{1:\infty}\|_1^m)/\|\Theta\|_p^{m-\a} \big]\\
&=&\E\big[ (\|1+\Theta_{1:\infty}\|_1^m-\|\Theta_{1:\infty}\|_1^m)/\|\Theta\|_p^{m-\a} \big]\,.
\eeao
Using the relation 
\beao
 \|1+\Theta_{1:\infty}\|_1^m-\|\Theta_{1:\infty}\|_1^m=\sum_{j=0}^{m-1} {m\choose j} \|\Theta_{1:\infty}\|_1^j\,,
\eeao
it is immediate that $\E\big[|W^{(p)}|^m]<\infty$ if we can show that
\beam\label{condition:m-moment:timechange} 
\E\big[ \|\Theta_{0:\infty}\|_1^{m-1}/\|\Theta\|_p^{m-\a}  \big]<\infty\,.
\eeam

\section{Application to \regvary\ stationary \ts\ models}\label{sec:exam}\setcounter{equation}{0}
\subsection{Models with \asy\ independence}
Such \regvary\ models are characterized by the fact that $\bfTh_t=\bf0$ for all 
$t\ne 0$. In this case, the same limit theory applies as for an iid \seq\ 
with the same marginal \ds\ as $\bfX$ provided the mixing and anti-clustering conditions of  Theorem~\ref{pr:prophybridch} are satisfied. A simple example is the following.
 
\bexam {\bf A stochastic volatility model.} \rm We consider the real-valued 
stationary process $X_t=\sigma_t\,Z_t$ where the \regvary\ iid \seq\ $(Z_t)$
is independent of the stationary \seq\ $(\sigma_t)$ of positive \rv s. The \seq\
$(X_t)$ is \regvary\ with index $\a>0$ and spectral tail process $\Theta_t=0$ for
$t\ne 0$. The anti-clustering and mixing conditions of Theorem~\ref{pr:prophybridch}
can be verified under mild conditions on the dependence structure of $(\sigma_t)$;
see Chapters 9 and 10 of \cite{mikosch:wintenberger:2024}.
\eexam

\subsection{Linear processes}\label{sec:linear}
Regularly varying linear processes have the rare property that their spectral 
tail processes are ``almost deterministic''. We will use this property 
in combination with the results in the iid \regvary\ case obtained by Logan et al. \cite{logan:mallows:rice:shepp:1973}.
\par
Consider the univariate causal linear process $(X_t)$ given by 
\beam\label{def:linear:ps}
 X_t = \sum_{j=0}^\infty \psi_j\, Z_{t-j}\,,\qquad t\in \Z\,,
\eeam
where $(Z_t)$ is iid \regvary\ with index $\alpha>0$ and
tail balance coefficients $q_{\pm}$, and $(\psi_j)$ is a real-valued \seq\ satisfying the condition, for an arbitrarily small $\delta>0$, 
\beao
\sum_{j=0}^\infty |\psi_j|^{(\a-\delta)\wedge 1}<\infty\,.
\eeao
For convenience we write $\psi_j=0$ for $j<0$.
If $\E[|Z|]<\infty$ we also assume $\E[Z]=0$. Under these conditions
the infinite series \eqref{def:linear:ps} converges a.s.; see 
Proposition 5.2.1 in \cite{mikosch:wintenberger:2024}. 
Section 5.2.6.1 in  \cite{mikosch:wintenberger:2024} yields the spectral tail process
\[
 \Theta_t:= \Theta_Z \dfrac{\psi_{J+t}}{|\psi_J|},\qquad t \in\bbz\,,
\]
where $\Theta_Z$ and $J$ are independent with distributions 
\[
 \P(\Theta_Z =\pm 1) =q_{\pm},\qquad \P(J=j) = \dfrac{|\psi_j|^\alpha}{\|\psi\|_\alpha^\alpha},\quad j\in\Z.
\]
Then
\[
 \|\Theta\|_\alpha^\alpha =\sum_{t\in\Z} |\Theta_t|^\alpha = \sum_{t\in \Z} \dfrac{|\psi_{J+t}|^\alpha}{|\psi_J|^\alpha} = 
\dfrac{\|\psi\|_\alpha^\alpha}{|\psi_J|^\alpha}, 
\]
the spectral cluster process $Q=\Theta/\|\Theta\|_\alpha$ is given by
\[
 Q_t = \dfrac{\Theta_t}{\|\Theta\|_\alpha} = \Theta_Z\dfrac{\psi_{J+t}}{\|\psi\|_\alpha},\qquad t \in \Z , 
\]
and for $p>\a$,
\beao
 \sum_{t\in \Z} Q_t &=& \Theta_Z \frac{\sum_{t\in \Z} \psi_t}{\|\psi\|_\alpha},
\quad \|Q\|_p^p = \sum_{t\in \Z}|Q_t|^p = \frac{\|\psi\|^p_p}{\|\psi\|_\alpha^p}\,,\\
\frac{\sum_{t\in \Z}Q_t}{\|Q\|_p}&=&\Theta_Z \frac{\sum_{t\in \Z}\psi_t}{\|\psi\|_p}\,,\quad
\max_{t\in\bbz}|Q_t|=\frac{\max_{j\ge 0}|\psi_j|}{\|\psi\|_\alpha},\qquad t \in \Z .
\eeao
Therefore the joint \chdf-Laplace 
transform $\Psi_X$ of \\ $(\xi_\alpha^X,\zeta_{\a,\infty}^X, (\zeta_{\alpha,p}^X)^p)$  for a  
linear causal process $(X_t)$ and $p>\a$ is given
 by
\beao
\Psi_X(u,x,\la)&=&
\exp \Big(\int_0^\infty \E[\ex^{iu y \frac{\sum_{j\in \Z}\psi_j}{\|\psi\|_\alpha} \Theta_Z
-\lambda y^p \frac{\|\psi\|_p^p}{\|\psi\|_\alpha^p} }\,\1\big(y\,\max_{j\ge 0}|\psi_j|/\|\psi\|_\a\le x\big)\\&&-1-i\,u\,y \frac{\sum_{j\in\Z}\psi_j}{\|\psi\|_\alpha}\Theta_Z {\bf 1}_{(1,2)}(\alpha) \big] d(-y^{-\alpha}) \Big). 
\eeao
We mention that the mixing and anti-clustering conditions for this result are automatically satisfied in the cases where $\psi_j=0$, $|j|>m$. Then Theorem \ref{pr:prophybridch} applies to these moving average processes for every fixed $m\ge 1$. The joint convergence result of Theorem \ref{pr:prophybridch}  
extends by passing $m\to \infty$. This is achieved by  an application of Billingsley's argument; see Theorem 3.2 of \cite{billingsley:1999}, and Section 9.3.2 of \cite{mikosch:wintenberger:2024} for details.

Writing $\Psi_Z$ for the corresponding \fct\ in the case of an iid \regvary\ \seq\ $(Z_t)$, we see that we have the relation
\beam\label{eq:readoff}
\Psi_X(u,x,\la)= \Psi_Z\Big(u\,\dfrac{\sum_{j\in \Z}\psi_j}{\|\psi\|_\alpha}\,,x\,
\dfrac{\|\psi\|_\a}{\max_{j\ge 0}|\psi_j|}\,,\la\,\dfrac{\|\psi\|_p^p}{\|\psi\|_\alpha^p}\,
\Big)\,.
\eeam 
If we write $(\xi_\alpha^Z,\zeta_{\a,\infty}^Z, (\zeta_{\alpha,p}^Z)^p)$ for the corresponding
quantities of the iid $Z$-process we can read off the following identity in 
law from \eqref{eq:readoff}:
\beao
\big(\xi_\alpha^X,\zeta_{\a,\infty}^X,\zeta_{\alpha,p}^X\big)
\eqd \Big(\dfrac{\sum_{j\in \Z}\psi_j}{\|\psi\|_\alpha}\,\xi_\a^Z\,,\dfrac{\max_{j\ge 0}|\psi_j|}{\|\psi\|_\a}\zeta_{\a,\infty}^Z\,,
\dfrac{\|\psi\|_p}{\|\psi\|_\alpha}\,\zeta_{\a,p}^Z
 \Big)\,.
\eeao
In words, the right-hand quantities are just re-scaled versions 
of the limit variables in the case of an iid \regvary\ \seq\ $(Z_t)$.
This observation is in agreement with work of Davis and Resnick \cite{davis:resnick:1985,davis:resnick:1985a,davis:resnick:1986} who proved limit theory 
for \fct als acting on linear processes with \regvary\ noise $(Z_t)$.
The main tools in the aforementioned papers are continuous mapping techniques
acting on converging \pp es.
\par
An immediate con\seq\  is the following ratio limit relations.

\bco\label{prop:linear:density}
Consider a \regvary\ causal linear process with index $\a\in (0,2)\backslash\{1\}$
satisfying the aforementioned conditions on $Z$ and $(\psi_j)$. Then for every $ p\in(\a,\infty]$,
\beao
\dfrac{S_n}{\gamma_{n,p}} 
\std  
\dfrac{\sum_{j\in \Z}\psi_j}{\|\psi\|_p}\dfrac{\xi_\a^Z}{\zeta_{\a,p}^Z}\,, \qquad \nto\,,\eeao
and the limit ratio has a moment generating \fct\ on the whole real line.
We consider some special cases.
\begin{enumerate}
\item[\rm 1.]
If $p=1$, $\a\in (0,1)$ and  $q_+q_->0$ then the limit ratio has a continuous 
density on $\|\psi\|_1^{-1}(-\sum_{j\in \Z}\psi_j,\sum_{j\in \Z}\psi_j)$ and singularities at the boundaries.
\item[\rm 2.]
If $p=2$ and $\alpha\in (1,2)$ then  
the limit ratio has a continuous density 
whose tails are \asy ally equivalent
to a Gaussian density.
\end{enumerate}
\eco
\begin{proof} By definition the limit ratio
\beao
\dfrac{\xi_\a^Z}{\zeta_{\a,p}^Z}=\dfrac{\xi_\a^Z}{\big\|\big(\Gamma_i^{-1/\a}\big)_{i\ge 1}\|_p}\le \dfrac{\xi_\a^Z}{\big\|\big(\Gamma_i^{-1/\a}\big)_{i\ge 1}\|_\infty}\,,
\eeao
and the existence of a moment-generating function follows from Lemma~\ref{lem:1} applied to the case $p=\infty$.
The comments about the limit densities in parts 1 and 2 follow from
the results in Logan et al. \cite{logan:mallows:rice:shepp:1973} mentioned in Section~\ref{subsec:iid}. \end{proof}

\begin{figure}[thbp]
\centerline{
\epsfig{figure=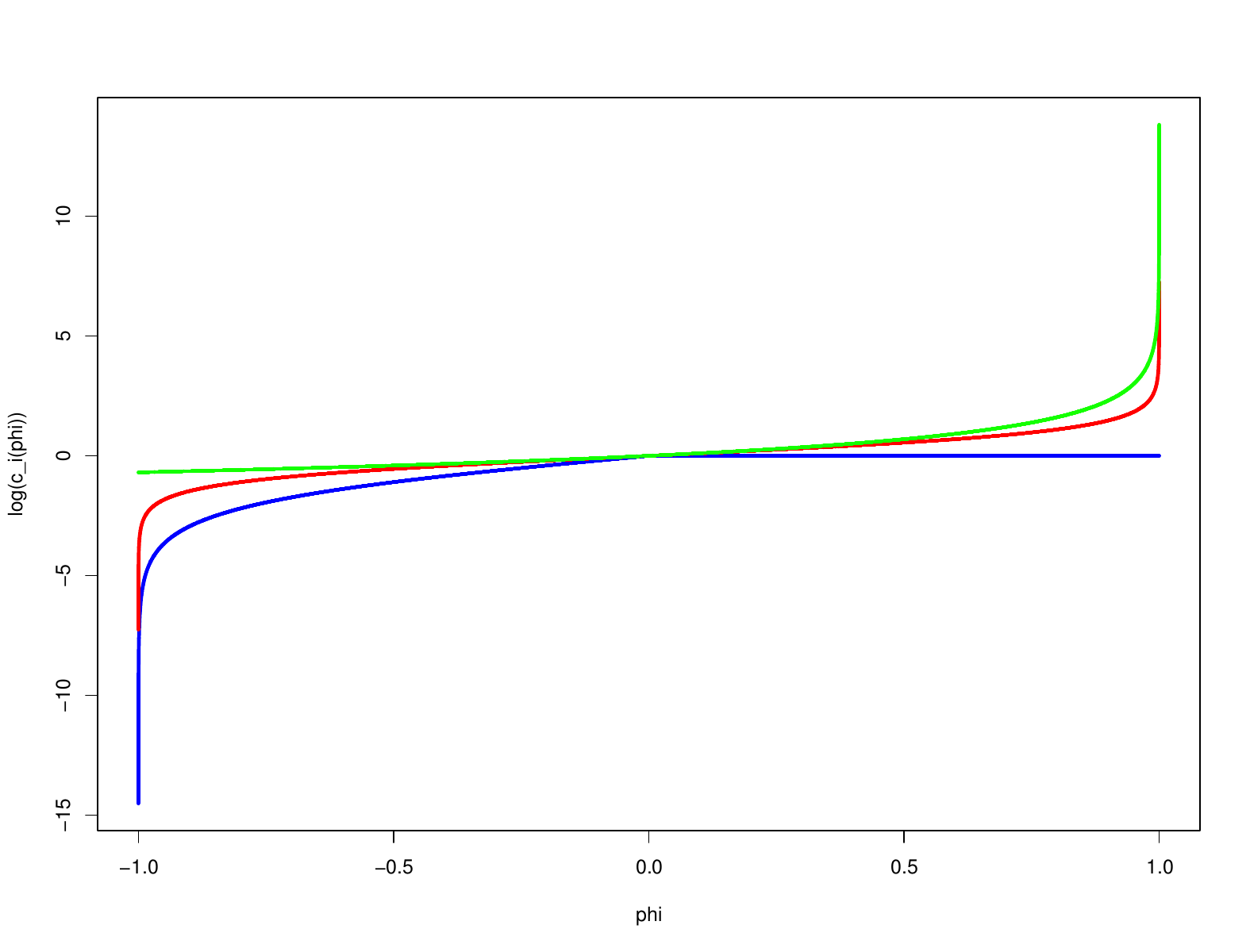,height=10cm,width=7cm,angle=-90}
}
\caption{The \fct s $c_1(\varphi)$ (blue), $c_2(\varphi)$ (red) and $c_\infty(\varphi)$ (green) on log-scale.}\label{fig:c2phi}
\end{figure}

\bexam {\bf A \regvary\ AR(1) process.}
\rm
Consider a \regvary\ stationary AR(1) process $(X_t)$ given as the causal 
solution to the difference equation
$X_t=\varphi X_{t-1}+Z_t$, $t\in \Z$, where $|\varphi|<1$ and 
$(Z_t)$ is iid \regvary\ with index $\alpha\in (1,2)$. Then we have the 
linear process \rep\ $X_t=\sum_{j=0}^\infty \varphi^j\,Z_{t-j}$, $t\in\bbz$. 
In particular,
\begin{align*}
& \dfrac{\sum_{j\in\bbz}\psi_j}{\|\psi\|_1}=\dfrac{1-|\varphi|}{1-\varphi}=:c_1(\varphi)\,,\quad 
\dfrac{\sum_{j\in \Z}\psi_j}{\|\psi\|_2}=
 \dfrac{\sqrt{1+\varphi}}{\sqrt{1-\varphi}}=:c_2(\varphi)\,, \\
&\quad
 \dfrac{\sum_{j\in\bbz}\psi_j}{ \max_{j\ge 0}|\psi_j|}=\dfrac 1 {1-\varphi}=:c_\infty(\varphi)\,.
\end{align*}
In Figure~\ref{fig:c2phi} we show the graphs of $\log (c_i(\varphi))$, $i=1,2,\infty$, indicating the magnitude of the deviation from the iid case.

\begin{figure}[htbp]
\centerline{
\epsfig{figure=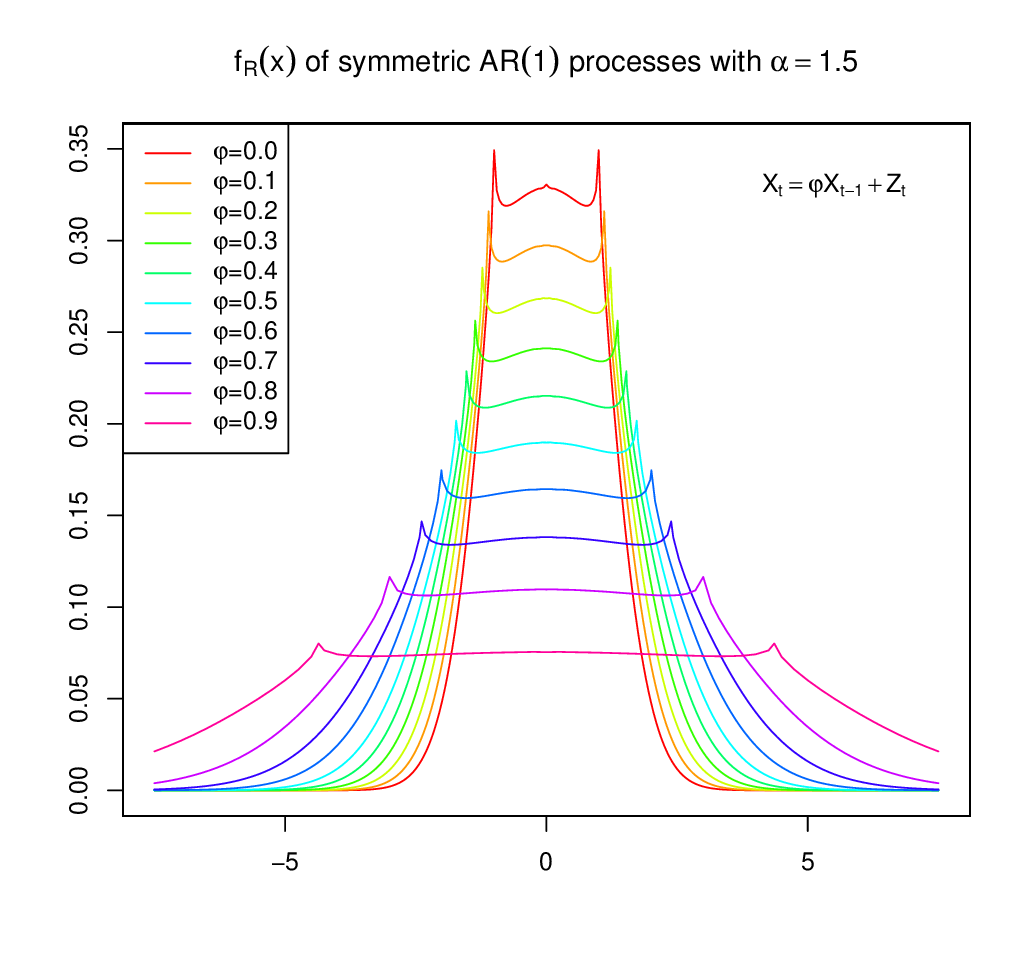,height=10cm,width=11cm}}
 \caption{Limit density of the studentized sums of a \regvary\
AR(1) process $X_t=\varphi\,X_{t-1}+Z_t$, $t\in\bbz$,  
for different values of $\varphi$ and symmetric \regvary\ $Z$ with index $\a=1.5$.}\label{fig:graph:R:ar1}
\end{figure}

Write $f_{R_{\a,2}^Z}$ and $f_{R_{\a,2}^X}$ for the respective densities of the limit ratios
$R_{\a,2}^Z=\xi_\a^Z/\zeta_{\a,2}^Z$ and 
$R_{\a,2}^X=\xi_\a^X/\zeta_{\a,2}^X \eqd  c_2(\varphi) R_{\a,2}^Z$
of the studentized \seq s 
$\sum_{t=1}^n Z_t/\big(\sum_{t=1}^n Z_t^2\big)^{1/2}$, $n\ge 1$, and 
$S_n/\gamma_{n,2}$, $n\ge 1$.
In Figure \ref{fig:graph:R:ar1} we draw the limit densities   
$f_{R_{\a,2}^X}(y) = |c_2(\varphi)|^{-1} f_{R_{\a,2}^Z}(y/c_2(\varphi))$ 
for an iid symmetric $Z$-\seq\ with $\alpha=1.5$ and various $\varphi$-values.
For the evaluation of $f_{R_{\a,2}^Z}$ we used the numerical method described in Appendix~\ref{calc:density:R:iid}. 
\eexam

\subsection{The solution to an affine \sre } We consider the causal solution to the 
affine \sre\ $X_t=A_t\,X_{t-1}+B_t$, $t\in\bbz$, where 
$(A_t,B_t)$, $t\in\bbz$, is an $\bbr^2$-valued iid \seq\ with possible dependence
between $A_t$ and $B_t$. Under general conditions on the \ds\ of a generic element
$(A,B)$ the solution to this equation exists and is \regvary\ with positive index; see Buraczewski et al. \cite{buraczewski:damek:mikosch:2016}, Mikosch and 
Wintenberger \cite{mikosch:wintenberger:2024}.
\begin{itemize}
\item {\bf Kesten-Goldie conditions} \cite{kesten:1973,goldie:1991}. Assume that there exists  a
positive solution $\a>0$ to the equation $\E[|A|^\a]=1$,
$\E[|A|^\a\log^+ |A|]<\infty$, $\E[|B|^\a]<\infty$, the conditional law of $\log |A|$
given $\{A\ne 0\}$ is non-lattice and $\P(A\,x+B)<1$ for every $x\in \bbr$.
Then $(X_t)$ is \regvary\ with index $\a$. In particular, $\P(\pm X_0>x)\sim c_{\pm} x^{-\a}$ as 
$\xto$ for constants $c_{\pm}$ \st\ $c_++c_->0$, and $c_+=c_-$ if $\P(A<0)>0$.  
The forward process is given by $\P(\Theta_0=\pm 1)= c_{\pm}/(c_++c_-)$, 
$\Theta_0$ is independent of $(A_t)$, and   
\beam\label{eq:Thetat}
\Theta_t=\Theta_0\,A_1\cdots A_t\,,\qquad t\ge 0.
\eeam
\item {\bf Grincevi\v cius-Grey conditions} \cite{grincevicius:1975,grey:1994}. Assume that $(B_t)$ is 
\regvary\ with index $\a>0$, $A\ge 0$ a.s., $\P(A=0)<1$, $\E[A^\a]<1$ and 
$\E[A^{\a+\delta}]<\infty$ for some $\delta>0$. Then 
$(X_t)$ is \regvary\ with index $\a$. In particular,  
$\P(\pm X_0>x)\sim (1-\E[A^\a])^{-1}\,\P(\pm B>x)$ as $\xto$. The 
forward process is given by $\P(\Theta_0=\pm 1)=q_{\pm}$ where
$q_{\pm}$ are the tail-balance coefficients of $B$, $\Theta_0$ is independent of
$(A_t)$, and $(\Theta_t)_{ t\ge 1}$ is given by  \eqref{eq:Thetat}.
\end{itemize}
\par
In what follows, we study the existence of moments of $R_{\a,p}$ for the solution of an affine \sre\
under the conditions of Kesten-Goldie or Grincevi\v cius-Grey. To achieve this we will apply
Theorem \ref{thm:momentsratio}. In particular, we will verify condition \eqref{condition:m-moment:timechange}; we recall it for convenience of the reader: 
\beam\label{condition:m-moment:timechangea} 
\E\big[ \|\Theta_{0:\infty}\|_1^{m-1}/\|\Theta\|_p^{m-\a}  \big]<\infty\,.
\eeam

\bpr\label{prop:sre:moment}
Consider the stationary solution $(X_t)$ to the affine \sre\ $X_t=A_t\,X_{t-1}+B_t$, $t\in\bbz$, and assume one of the following 
conditions.
\begin{enumerate}
\item[\rm 1.] The conditions of the Kesten-Goldie theory,
in particular,  $\E[|A|^\a]=1$ for some $\a\in (1,2)$.
\item[\rm 2.] 
The conditions of the Grincevi\v cius-Grey theory. In particular, $B$ is \regvary\ with index $\a\in (1,2)$, 
$A\ge 0$ a.s. and $\E[A^\a]<1$. In addition, we require one of the following conditions:\\
(i) The \sre\ $Z_t=A_t\,Z_{t-1}+1$, $t\in\bbz$,
satisfies the Kesten-Goldie conditions for some $\beta>\a$, in particular, $\E[A^\beta]=1$.\\
(ii) $\E[A^\beta]<1$ for every $\beta>0$.
\end{enumerate}
Then $(X_t)$ is \regvary\ with index $\a$ and all moments $\E\big[R_{\a,p}^m\big]$, $m>0$, $p>\a$, are finite.
\epr

\bre\label{remark:condition:m-moment} The case $m=2$ is immediate by verifying \eqref{condition:m-moment:timechangea}.
We have $\Theta_t= \Theta_0 A_1\cdots A_t$, $t\ge 1$, $|\Theta_0|=1$, and 
\beao
\E\big[ \|\Theta_{0:\infty}\|_1 /\|\Theta\|_p^{2-\a}  \big]&\le & 
\E\big[ \|\Theta_{0:\infty}\|_1]=\sum_{t\ge 0} (\E[|A|])^t<\infty\,.
\eeao 
The property $\E[|A|]<1$ follows from convexity of
$f(h)=\E[|A|^h]$, $h>0$, and  $ \E[|A|^\a]\le 1$ for $\a\in (1,2)$. 
\ere
\begin{proof} {\em Part 1.}
The \regvar\ of $(X_t)$ follows from the Kesten-Goldie theory. The absolute values of the forward spectral tail process 
are given by $|\Theta_t|=|A_1\cdots A_t|$ for $t\ge 1$ and $|\Theta_0|=1$ a.s. 
In view of Theorem~\ref{thm:momentsratio}  and the aforementioned arguments it suffices to show that \eqref{condition:m-moment:timechangea} holds for every
even positive integer $m$. It is implied by
\beao
 \E[\|\Theta_{0:\infty}\|_1^{m-1}/\|\Theta_{0:\infty}\|_\infty^{m-\alpha}]<\infty\,.
\eeao
For a random variable $A\stackrel{d}{=}A_1$ which is independent of $(\Theta_t)_{t\ge 0}$, we have 
\beao
\big(\|\Theta_{0,\infty}\|_1, \|\Theta_{0,\infty}\|_\infty\big) &=& 
\Big(1+|A_1|\Big(1+\sum_{t=2}^\infty |A_2\cdots A_t|\Big)\,,1\ \vee |A_1| \Big(1\vee \max_{t\ge 2}\,\big(|A_2\cdots A_t| \Big)\Big) \\
 &\stackrel{d}{=} & (1+|A|\,\|\Theta_{0:\infty}\|_1\,, 1 \vee |A|\,\|\Theta_{0:\infty}\|_\infty).
\eeao
Therefore 
\begin{align*}
\frac{\|\Theta_{0:\infty}\|_1^{m-1}}{\|\Theta_{0:\infty}\|^{m-\alpha}_\infty} & \stackrel{d}{=} 
\frac{(|A| \|\Theta_{0:\infty}\|_1+1)^{m-1}}{(1\vee |A|\|\Theta_{0:\infty}\|_\infty)^{m-\alpha}} \\
&= \Big(
\frac{|A| \|\Theta_{0:\infty}\|_1+1 }{ (1\vee |A|\|\Theta_{0:\infty}\|_\infty)^{(m-\alpha)/(m-1)}}
\Big)^{m-1} \\
&\le \Big(
\frac{|A| \|\Theta_{0:\infty}\|_1}{(|A|\|\Theta_{0:\infty}\|_\infty)^{(m-\alpha)/(m-1)}}+1
\Big)^{m-1} \\
& = \Big(
|A|^{(\alpha-1)/(m-1)} \Big(
\frac{\|\Theta_{0:\infty}\|_1^{m-1}}{\|\Theta_{0:\infty}\|^{m-\alpha}_\infty}
\Big)^{1/(m-1)}+1 
\Big)^{m-1}, 
\end{align*}
and the random variable $(\|\Theta_{0:\infty}\|_1^{m-1}/\|\Theta_{0:\infty}\|^{m-\alpha}_\infty)^{1/(m-1)}$ is 
stochastically dominated  by the marginal \ds\ of the solution to the
affine \sre\ 
\begin{align}
\label{SRE:Y:1/m-1}
Y_t^{1/(m-1)}= |A_t|^{(\a-1)/(m-1)}\,Y_{t-1}^{1/(m-1)}+1\,,\qquad  t\in\bbz\,. 
\end{align}
By the Kesten-Goldie theory it has a \regvary\ stationary solution
$(Y_t^{1/(m-1)})$ satisfying the identity in law $Y^{1/(m-1)}\eqd  |A|^{(\a-1)/(m-1)}\,Y^{1/(m-1)}+1$ where a generic element $Y^{1/(m-1)}$ and $A$ are independent.
But $\E\big[\big(|A|^{(\a-1)/(m-1)}\big)^\beta\big]=1$ has the unique solution $\beta=\a(m-1)/(\a-1)$. Therefore $Y$ is \regvary\ with index $\a/(\a-1)>1$
and has a finite first moment. Stochastic domination implies that $\E[\|\Theta_{0:\infty}\|_1^{m-1}/\|\Theta_{0:\infty}\|_\infty^{m-\a}]<\infty$ as well and the desired result follows for every $p>\a$.\\[1mm]
{\em Part 2.} The \regvar\ of $(X_t)$ with index $\a\in (1,2)$ follows from the Grincevi\v cius-Grey theory. The forward spectral
tail process is given by $\Theta_t=A_1\cdots A_t$, $t\ge 1$, and $|\Theta_0|=1$ a.s.
\par
In the 
the case $(i)$ when $\E[A^\beta]=1$ for some $\beta< m$,
$\alpha\in (1,2)$ and the conditions of Kesten-Goldie are satisfied for $(A_t)$ (then, in particular, $\E[A^\a]<1$ and $\beta>\a$)
and $B_t=1$, $t\in\bbz$, then the affine \sre\ $Z_t=A_t\,Z_{t-1}+1$, $t\in\bbz$, has a \regvary\ stationary solution with index $\beta$. Proceeding as in the proof of Part 1 (with $\a$ replaced by $\beta$), we may conclude that
the affine \sre\ $Y_t^{1/(m-1)}=A_t^{(\beta-1)/(m-1)}\, Y_{t-1}^{1/(m-1)}+1$ has a \regvary\ stationary solution $(Y_t)$ with index $\beta/(\beta-1)$, hence it has the first moment and 
\beao
\infty>\E\big[\|\Theta_{0:\infty}\|_1^{m-1}/\|\Theta_{0:\infty}\|_\infty^{m-\beta}\big]\ge \E\big[\|\Theta_{0:\infty}\|_1^{m-1}/\|\Theta_{0:\infty}\|_\infty^{m-\a}\big]\,,
\eeao
Therefore \eqref{condition:m-moment:timechangea} is satisfied.

The case $(ii)$ is handled by a similar argument as in Remark~\ref{remark:condition:m-moment} 
with general even integers~$m$.
\end{proof}

\subsection{The \garch\ model}
The Generalized Autoregressive Conditionally Hete\-ro\-sce\-das\-tic\ (GARCH)  process is one of the standard models for returns of prices in financial  \tsa ; see Engle \cite{engle:1982}, 
Bollerslev \cite{bollerslev:1986}. Here we focus on the \garch\ model introduced in \cite{bollerslev:1986}: $X_t=\sigma_t\,Z_t$, $t\in\bbz$, where $(Z_t)$ is an iid \seq\
whose generic element $Z$ satisfies $\E[Z]=0$, $\E[Z^2]=1$. The {\em squared volatility process} $(\sigma_t^2)$ is the stationary solution to the \sre\
$\sigma_t^2=\a_0+(\a_1\,Z_{t-1}^2+\beta_1)\,\sigma_{t-1}^2=:B_t+A_t\,\sigma_{t-1}^2$, $t\in\bbz$, where $\a_0>0$ and $\a_1,\beta_1$ are non-negative parameters satisfying 
$\E[\log (\a_1\,Z^2+\beta_1)]<0$. If the \sre\ defining $(\sigma^2_t)$ satisfies the conditions of the Kesten-Goldie theory, in particular
$\E[(\a_1\,Z^2+\beta_1)^\a] =1$ for a positive $\a$, then the squared volatility process is \regvary\ with index $\a$, and $(X_t)$ inherits \regvar\ with 
index $2\,\a$ and forward spectral tail process $(\Theta_t)$ given under a change of \ms\ by
\beao
\P\big( (\Theta_0,\ldots,\Theta_h)\in \cdot \big) =\E\Big[\dfrac{|Z_0|^{2\a}}{\E[|Z_0|^{2\a}]}\,\1\big(|Z_0|^{-1}\,
\big(Z_0\,,\Pi_1^{1/2}\,Z_1,\ldots,\Pi_h^{1/2}\,Z_h\big)\in \cdot \big)\Big]\,, \quad h\ge 0\,,
\eeao
where $\Pi_t=A_1\cdots A_t$, $t\ge 1$, $\Pi_0=1$.
We refer to Section 5.6.3.6 in \cite{mikosch:wintenberger:2024} for these results and the conditions under which they hold, for example they are satisfied for 
standard normal $Z$.
\bpr Consider a \garch\ process as described above.
Assume the conditions of the Kesten-Goldie theory for $(\sigma_t^2)$. In particular, $\E[(\a_1\,Z^2+\beta_1)^\a] =1$  for some $\a\in (0.5,1)$. Then $(X_t)$ is \regvary\
with index $2\,\a\in (1,2)$  and all moments $\E[R_{2\a,p}^m]$, $m>0$, $p>2\,\a$, are finite.
\epr 
\begin{proof}
Assuming $2\a\in (1,2)$ and that $m\ge 2$ is an even integer,
we will verify the condition \eqref{condition:m-moment:timechangea} for $p>2\a$. It implies $\E[R_{2\,\a,p}^m]< \infty$.  We have
\beao
\E\Big[\dfrac{\|\Theta_{0:\infty}\|_1^{m-1}}{\|\Theta_{0:\infty}\|_p^{m-\red 2\a}}\Big]
&=&\E\Big[\dfrac{|Z_0|}{\E[|Z_0|^{2\a}]}\,\Big(
\dfrac{\sum_{t=0}^\infty |Z_t|\,\Pi_t^{1/2}}{\Big( \sum_{t=0}^\infty |Z_t|^p\,\Pi_t^{p/2}\Big)^{(m-\red 2\a)/(p(m-1))}}\Big)^{m-1}\Big]\\
&\le &
\E\Big[\dfrac{|Z_0|}{\E[|Z_0|^{2\a}]}\,\Big(
\dfrac{\sum_{t=0}^\infty |Z_t|\,\Pi_t^{1/2}}{ \max_{t\ge 0} (|Z_t|\,\Pi_t^{1/2})^{(m-\red 2\a)/(m-1)}}\Big)^{m-1}\Big]\,.
\eeao
Write
$\Pi_{s,t}=A_s\cdots A_t$, $t\ge s$, $\Pi_{s,t}=1$ if $t<s$.
We observe that
\beao
 \Big(
\sum_{t=0}^\infty |Z_t| \Pi_t^{1/2},\,\max_{t \ge 0} \big(|Z_t| \Pi_t^{1/2} 
\big) \Big) = \Big( A_1^{1/2}
\sum_{t=1}^\infty |Z_t|\Pi_{2,t}^{1/2} +|Z_0|,\, |Z_0| \vee A_1^{1/2} \max_{t\ge 1} \big(|Z_t| \Pi_{2,t}^{1/2}\big)
\Big). 
\eeao
In what follows, $c$ stands for positive constants whose values are not of interest.
Therefore for some constant $c>0$,
\begin{align}
& |Z_0|\Big(
\frac{\sum_{t=0}^\infty |Z_t| \Pi_t^{1/2}}{\max_{t \ge 0} \big(|Z_t| \Pi_t^{1/2} 
\big)^{(m-{\red 2\alpha})/(m-1)}}
\Big)^{m-1} \nonumber \\
& \le |Z_0|\Big(
A_1^{0.5 (\red 2\alpha-1)/(m-1)} \frac{\sum_{t=1}^\infty |Z_t| \Pi_{2,t}^{1/2} }{\max_{t \ge 1} \big(|Z_t| 
\Pi_{2,t}^{1/2}
\big)^{(m-{\red 2\alpha})/(m-1)}} +|Z_0|^{({\red 2\alpha}-1)/(m-1)} \Big)^{m-1} \nonumber\\
& \le  c\,(|Z_0|+|Z_0|^{\alpha})  \,\Big(A_1^{0.5 (\red 2\alpha-1)} \Big(
\dfrac{\sum_{t=1}^\infty |Z_t| \Pi_{2,t}^{1/2}}{\max_{t \ge 1} \big(|Z_t| 
 \Pi_{2,t}^{1/2}
\big)^{(m-{\red 2\alpha})/(m-1)}}
\Big)^{m-1} + 1\Big)\,. \nonumber 
\end{align}
Taking expectations, we find some constant $c>0$ \st 
\beao
 \E\Big[
\frac{\|\Theta_{0:\infty}\|_1^{m-1}}{\|\Theta_{0:\infty}\|_\infty^{m-{\red 2\alpha}}}
\Big] &\le& c\,\Big(
\E\Big[\Big(
\frac{\sum_{t=0}^\infty |Z_t| \Pi_{t}^{1/2}}{\max_{t \ge 0} \big(|Z_t| 
 \Pi_{t}^{1/2}
\big)^{(m-{\red 2\alpha})/(m-1)}}
\Big)^{m-1} \Big]+1\Big)\,.
\eeao
But we have
\beao\lefteqn{
 \frac{\sum_{t=0}^\infty |Z_t| \Pi_t^{1/2}}{\max_{t \ge 0} \big(|Z_t| \Pi_t^{1/2} 
\big)^{(m-{\red 2\alpha})/(m-1)}}}\\ &  \stackrel{d}{=}& \frac{A_0^{1/2} \sum_{t=0}^\infty |Z_t|\Pi_t^{1/2} +|Z_{-1}|}{
\big(
|Z_{-1}| \vee {A_0}^{1/2}\max_{t\ge 0} (|Z_t| \Pi_t^{1/2})
\big)^{(m-{\red 2\alpha})/(m-1)}} \\
& \le& A_{0}^{0.5\,({\red 2\alpha}-1)/(m-1)} \frac{\sum_{t=0}^\infty |Z_t| \Pi_t^{1/2}}{\max_{t \ge 0} \big(|Z_t| \Pi_t^{1/2} 
\big)^{(m-{\red 2\alpha})/(m-1)}} + |Z_{-1}|^{({\red 2\alpha}-1)/(m-1)}. 
\eeao
The quantity on the \lhs\ is stochastically dominated by 
the marginal distribution of the solution to the affine stochastic recurrence equation
\[
 Y_t^{1/(m-1)} = A_t^{0.5({\red 2\alpha}-1)/(m-1)} Y_{t-1}^{1/(m-1)} + C_t^{({\red 2\alpha}-1)/(m-1)},
\]
where $(A_t,C_t)$, $t\in\bbz$,  is an iid sequence such that 
$A_t=\alpha_1 Z_t^2 +\beta_1$ and $C_t=|Z_{t-1}|$. A generic element $(A,C)$ satisfies the equation
$\E[(A^{\beta ({\red 2\alpha}-1)/(m-1)})]=1$ with $\beta=\alpha (m-1)/({\red 2\alpha}-1)$ and $\E[(C^{2\beta ({\red 2\alpha}-1)/(m-1)})]<\infty$. 
By the Kesten-Goldie theory it has a regularly varying stationary solution $(Y_t^{1/(m-1)})$ satisfying the identity in law 
$Y^{1/(m-1)}\stackrel{d}{=} A^{({\red 2\alpha}-1)/(m-1)}Y^{1/(m-1)} + C^{({\red 2\alpha}-1)/(m-1)}$ where a generic element $Y^{1/(m-1)}$ and 
$(A,C)$ are independent. We note that $Y^{1/(m-1)}$ has finite moments of order $m-1<2\beta$ and we conclude with an argument similar to the proof of Proposition~\ref{prop:sre:moment}. 
\end{proof}

\subsection{An example with certain infinite moments}
In the previous sections we considered \regvary\ linear processes, solutions to affine \sre s and the \garch\ process. In these cases we were in the ideal situation that
all moments of the limit ratios $R_{\a,p}$, $p>\a$,  exist. Here we provide an example of a  \regvary\ stationary \seq\  where certain moments do not exist. We use the same
notation as in the previous sections.

\bpr\label{prop:nomoments}
There exists a real-valued \regvary\ stationary process $(X_t)$ with index $\a\in(0,2)\backslash\{1\}$, centered for $\a\in (1,2)$, with the following properties.
\begin{enumerate}
\item[\rm 1.]
The conditions of 
Theorem~\ref{pr:prophybridch} are satisfied. In particular,
for every $p>\a$, 
$S_n/\gamma_{n,p}\std \xi_\a/\zeta_{\a,p}=:R_{\a,p}$, $\nto$.
\item[\rm 2.] Assume $\a\in (0,1) $ and $p>1 $.   If the   integer $m$ satisfies
\beao
m>  \dfrac{2p-\a}{p-1}  \,,
\eeao
then we have $\E[|R_{\a,p}|^m]=\infty$.
\item[\rm 3.] Assume $\a\in (1,2) $ and $p> \a$. For any even integer $m$ satisfying 
\beao 
m> \Big(2\Big(\dfrac{\a-1}{2-\a}\vee 1\Big)p-\a\Big)\Big/(p-1) \,,
\eeao
we have $\E[|R_{\a,p}|^{m'}]=\infty$ for all $m'>m$.
\end{enumerate}
\epr
The remainder of this section is devoted to the proof of Proposition~\ref{prop:nomoments}.
\subsection*{Construction of $(X_t)$}
Define a positive integer-valued  regenerative  
time-homogeneous  \MC\ $(K_t)$ as follows: choose a positive integer-valued \rv\ $L$ and set
\beao
\left\{
\barr{ll}
K_t\eqd L\,, &\qquad \mbox{if $K_{t-1}= 1$,}\\
K_t=K_{t-1}-1\,, &\qquad \mbox{otherwise.}
\earr
\right.
\eeao
It has the atom $\{1\}$ by definition. If $\E[L]<\infty$ then the Markov chain $(K_t)$ has a stationary version and the \ds\ of a generic element $K$ satisfies
\beao
\P(K\ge k)=\P(K\ge k+1)+\P(L\ge k)\P(K=1), \qquad k\ge 1\,,
\eeao
hence $\P(K=k)=\P(L\ge k)\,\P(K=1)$. 
Summing over $k\ge 1$ on both sides of this equality, we get $\P(K=1)=1/\E[L]$. 
In what follows, we assume $ \E[L^3]<\infty$. Then we also have
\beao
\E[L] \,\E[K]&=& \E[L]\,\sum_{k=1}^\infty k\,\P(K=k) = \E[L]\,\P(K=1) \,\sum_{k=1}^\infty k\,\P(L\ge k)   
 = \E[(L+1)L/2],
\eeao
and a similar argument yields $\E[K^2]<\infty$.
\par
Now consider a \regvary\ iid \seq\ $(Z_t)$ with index $\a\in (0,2)\backslash\{1\}$ independent of $(K_t)$. We assume that a generic element $Z$ is positive.
We define the max-moving average of random size 
\beao
X_t=\max_{1\le j\le K_t}Z_{t-1+j}-\E\big[\max_{1\le j\le K}Z_{j}\big]\1_{(1,2)}(\a)\,,\qquad t\in\bbz\,.
\eeao
It constitutes a stationary process with generic element $X$.

\subsection*{Regular variation of $(X_t)$}
 Since $\E[K]<\infty$ we have $\P(X>x)\sim \E[K]\,\P(Z>x)$\,, $\xto$; see  Proposition~B.2.8,1. in~\cite{mikosch:wintenberger:2024}. Hence the marginal
\ds s of $(X_t)$ are \regvary\ with index $\a$. To show \regvar\ of $(X_t)$ it suffices to calculate the forward spectral tail process; see Basrak and Segers
\cite{basrak:segers:2009}; cf. Theorem 4.2.5 in \cite{mikosch:wintenberger:2024}.

\subsubsection*{The forward spectral tail process}
Our aim is to determine the forward spectral tail process $(\Theta_t)$ defined  via the limit relations for $h\ge 0$,
\beam\label{eq:xx1}
 \P(|X_0|^{-1}(X_0,\ldots,X_h)\in\cdot\, | |X_0|>x)\stw \P((\Theta_0,\ldots,\Theta_h)\in\cdot)\,,\qquad \xto\,.
\eeam
We apply the single big jump principle in a similar way as in Section 5.2.5 in  \cite{mikosch:wintenberger:2024} given that $K_0=k$. Therefore by an application of Lemma 5.2.6 in \cite{mikosch:wintenberger:2024} we have
\beao
\dfrac{\P(|X_0|>x |K_0=k)}{\P(Z>x)}\sim \dfrac{\P(\max_{0\le t\le k-1}Z_t>x)}{\P(Z>x)}\sim \dfrac{\P(\sum_{ t=0}^{k-1}Z_t>x)}{\P(Z>x)}\to k\,,\qquad \xto\,,
\eeao
where we used the subexponentiality of the \ds\ of $Z$ in the last step.
The single big jump manifests as
\beao
\lim_{\xto} \sum_{j=0}^{k-1}\P\Big(Z_j>x,\max_{0\le i\neq j\le k-1}Z_i<x\,\Big|\, |X_0|>x\,,K_0=k\Big)=1\,.
\eeao 
Moreover, given $K_0=k$ we have that $X_t$ is independent of $X_0$ for $t\ge k$.
First we assume $h\ge k$ and observe that for any continuity sets $A\subset \R^k$ and $B\in\R^{h-k+1}$
\beao
p(A\times B,x,k,h)&:=&\P(|X_0|^{-1}(X_0,\ldots,X_h)\in A\times B \,\mid\, |X_0|>x,K_0=k)\\
&=&\P(|X_0|^{-1}(X_0,\ldots,X_{k-1})\in A\,\mid\, |X_0|>x,K_0=k)\\
&&\times \P(|X_0|^{-1}(X_k,\ldots,X_h)\in B\,\mid\, |X_0|>x,K_0=k)\,,
\eeao
where by independence the second probability converges to $\vep_{\mathbf{0}}(\cdot)$ with $\mathbf{0}$ an $(h-k+1)$-dimensional zero vector.
Here $\vep_v$ stands for Dirac \ms\ at any point $v$.
For the first probability we have as $x\to\infty$
\beao
 p(\cdot,x,k,k-1)&:=&\P(|X_0|^{-1}(X_0,\ldots,X_{k-1})\in\cdot\, \mid\, |X_0|>x,K_0=k)\\
&=&\dfrac{\P(|X_0|^{-1}(X_0,\ldots,X_{k-1})\in\cdot\, , |X_0|>x,K_0=k)}{ \P(|X_0|>x,K_0=k)}\\
&\sim & \sum_{j=0}^{k-1}\P\Big(\dfrac{\big(\max_{0\le t\le k-1}Z_t,\max_{1\le t\le k-1}Z_t,\ldots, Z_{k-1}\big)}{\max_{0\le t\le k-1} Z_t}\in\cdot\, ,Z_j>x \,,
\\
&&\qquad \max_{0\le i\neq j\le k-1}Z_i<x\Big)
\Big/
\sum_{j=0}^{k-1}\P\big(Z_j>x \,,\max_{0\le i\neq j\le k-1}Z_i<x\big)\\
&\sim& \sum_{j=0}^{k-1}\P\Big(\dfrac{\big(\max_{-j\le t\le k-1-j}Z_t,\max_{1-j\le t\le k-1-j}Z_t,\ldots,Z_{k-1-j}\big)}{\max_{-j\le t\le k-1-j}Z_t}\in\cdot\,, Z_0>x\,,\\&&\qquad \textcolor{red}{\max_{-j\le i\neq 0\le k-1-j}Z_i<x}\Big)\Big/ (k\P(Z_0>x))\,,
\eeao
where the last display follows from the stationarity of $(Z_t)$. Recall the spectral tail process of the iid $Z$-\seq ,  $\Theta^Z_t=\1(t=0)$, $t\in\bbz$. Therefore, as $\xto$, {\red $\P(\max_{-j\le i\neq 0\le k-1-j}Z_i<x\mid Z_0>x)\to 1$, and}
\beao
p(\cdot,x,k, k-1)
&\sim&\dfrac1k\sum_{j=0}^{k-1}\P\big(Z_0^{-1}\big(\max_{-j\le t\le k-1-j}Z_t,\max_{1-j\le t\le k-1-j}Z_t,\ldots,Z_{k-1-j}\big)\in\cdot\, \mid Z_0>x\big)\\
&\sim&\dfrac1k \sum_{j=0}^{k-1}\P\big(\big(\1(-j\le 0\le k-1-j),\ldots,\1(k-1-j=0)\big)\in\cdot\big)\\
&=&\dfrac1k \sum_{j=0}^{k-1}\vep_{(\1(t\le j))_{0\le t\le k-1}}(\cdot)\,.
\eeao
Now together with the second probability in $p(\cdot,x,k,h)$, for $h\ge k$ 
\beao
\P(|X_0|^{-1}(X_0,\ldots,X_{h})\in\cdot\, \mid\, |X_0|>x,K_0=k)\sim\dfrac1k \sum_{j=0}^{k-1}\vep_{(\1(t\le j))_{0\le t\le  h}(\cdot)}\,,\qquad \xto\,.
\eeao
This asymptotic relation also holds for $h<k$ since an approximation similar to that for $p(\cdot,k,k-1)$ works. 
Thus 
\beao\lefteqn{
\P(|X_0|^{-1}(X_0,\ldots,X_{h})\in\cdot\, \mid\, |X_0|>x)}\\
&=&\dfrac{ \sum_{k=1}^\infty\P(|X_0|^{-1}(X_0,\ldots,X_{h})\in\cdot\, \mid \,|X_0|>x,K_0=k)\P(|X_0|>x\mid K_0=k)\,\P(K_0=k)
}{\sum_{k=1}^\infty\P(|X_0|>x\mid K_0=k)\,\P(K_0=k)
}
\\
&\sim& \dfrac{\sum_{k=1}^\infty k^{-1}\sum_{j=0}^{k-1}\vep_{(\1(t\le j))_{0\le t\le h}}(\cdot)k\,\P(K=k)}{\sum_{k=1}^\infty k\, \P(K=k)}\\
&\sim&\dfrac{\E[\sum_{j=0}^{K-1}\vep_{(\1(t\le j))_{0\le t\le h}}(\cdot)]}{\E[K]}\,,\qquad \xto\,.
\eeao
By definition of the forward spectral tail process the \rhs\ coincides with \\
$\P\big((\Theta_0,\ldots,\Theta_h)\in \cdot\big)$; see \eqref{eq:xx1}.
Since $|X|$ is \regvary\ and $(X_t)$ admits a forward spectral tail process  $(X_t)$ is a \regvary\ \seq .

\subsection*{The anti-clustering condition \eqref{cond:acac}} 
Assuming $r_n=o(a_n^ 2/n\wedge n)$,  we have
\beao
nr_n(a_n^{-1}\E\big[\max_{1\le j\le K}Z_{j}\big]\1_{(1,2)}(\a))^2=o(1)\,,\qquad \nto\,.
\eeao
Hence we may neglect the centering term in the expression of $X_t$. Then \eqref{cond:acac} will follow if
\beao
\lim_{l\to\infty}\limsup_{\nto}  n \,\sum_{t=l}^{r_n}\E\Big[\underbrace{\Big(a_n^{-1}\max_{t\le j\le t + K_t-1}Z_j \wedge  1\Big)\,\Big(  a_n^{-1}\max_{0\le j\le K_0-1}Z_j\wedge
1\Big)}_{ =:A_{n;0,t}}\Big]=0\,.
\eeao
Using the regeneration property of the Markov chain $(K_t)$, we have
\beao 
&& \sum_{t=l}^{r_n}\E\big[\1(K_0\le t)
A_{n;0,t}\big] \\
&&= \sum_{t=l}^{r_n}\E\Big[\1(K_0\le t)\Big(  a_n^{-1}\max_{0\le j\le K_0-1}Z_j\wedge
1\Big)\Big]\E\Big[\Big(a_n^{-1}\max_{t\le j\le t + K_t-1}Z_j \wedge  1\Big)\Big]\\
&&\le r_n\E\Big[\Big(  a_n^{-1}\max_{0\le j\le K-1}Z_j\wedge
1\Big)\Big]^2=\textcolor{red}{o(n^{-1})}\,,\qquad \nto\,.
\eeao
Here we applied Karamata's theorem (see Bingham et al. \cite{bingham:goldie:teugels:1987}) for $\a\in(0,1)$ while we exploited the property $r_n=o(a_n^2/n)$ for $\a\in(1,2)$.
Thus it suffices to show that
\beam\label{eq:9} 
\lim_{l\to\infty}\limsup_{\nto}  n \,\sum_{t=l}^{r_n}\E\big[\1(K_0> t)\,A_{n;0,t}
\big]=0\,.
\eeam
If $K_0>t$ we have $K_t+t =K_0$. Therefore
\beao
\lefteqn{ \E\big[\1(K_0> t)\,A_{n;0,t}
\big]}\\
&\le& \E\Big[\1(K_0> t) \Big(  a_n^{-1}\max_{0\le j\le K_0-1}Z_j\wedge
1\Big)^2\Big]
\le \E\Big[\1(K_0> t)\sum_{j=0}^{K_0-1}(  a_n^{-1}Z_j\wedge
1)^2\Big]\\
&=& \E[\1(K> t)\,K]\E[(  a_n^{-1}Z\wedge
1)^2]
\le  c\,\E[\1(K> t)\,K]\P(Z>a_n)
\eeao
for some constant $c>0$ where we applied Karamata's theorem and used the independence between $(K_t)$ and $(Z_t)$.  For sufficiently large $n$ we have
\beao
n\,\sum_{t=l}^{r_n}\E\big[\1(K_0> t)\,A_{n;0,t}
\big] &\le&c\,\sum_{t=l}^{r_n}\E[\1(K> t)\,K]
\le c\,\E[(K-l)_+\,K]\,.
\eeao
Then \eqref{eq:9} follows from $\E[K^2]<\infty$ and \eqref{cond:acac} is satisfied for $r_n=o(a_n^ 2/n\wedge n)$.
\subsection*{The mixing condition \eqref{eq:wdepms2}}
The strong mixing coefficients $\a_t$ of $(X_t)$ satisfy the relation
$\a_t\le \P(K\ge t)$.
This follows from the independence of $(Z_t)$ and the fact that 
the $\sigma$-fields $\sigma(X_t, t\le 0)\subset \sigma(Z_t, K_t, t\le K_0-1)$ are  independent of $(X_t) _{ t\ge K_0}$ given $K_0$ due to the renewal structure of $(K_t)$. Therefore the mixing property \eqref{eq:wdepms2} holds for any integer sub\seq\
$(r_n)$ given there exists $(\ell_n)$ \st\ $\ell_n=o(r_n)$ and $(n/r_n) \a_{\ell_n}\to 0$ hold;  see Lemma 9.1.3 of \cite{mikosch:wintenberger:2024}.

\subsection*{Moment properties of the limit ratios} 
 Recall the aforementioned conditions on the intermediate 
\seq s $(r_n)$ and $(\ell_n)$ required above: 
\beao
r_n=o(a_n^2/n\wedge n)\,,\qquad \ell_n=o(r_n)\,,\qquad (n/r_n) \alpha_{\ell_n}\le (n/r_n)\P(K>\ell_n)
\to 0\,.
\eeao
Now we also assume that $K$ has all moments of order less than $q$ for some 
$q>2$ but $\E[K^q]=\infty$. Also recall that $a_n=n^{1/\a}l(n)$ for
a \slvary\ \fct\ $l(n)$. Taking into account all these requirements,
we get the existence of intermediate sequences $(r_n)$ and $(\ell_n)$ 
under the sufficient condition $q>2(\a-1)/(2-\a)$; it is restrictive 
only for $\a>1.5$. Now we can apply Theorems~\ref{pr:prophybridch} and 
\ref{thm:momentsratio}: 
$S_n/\gamma_{n,p}\std R_{\a,p}$ for $p>\a$, and the moment properties of $R_{\a,p}$ depend on those of  $W^{(p)}=\sum_{t\in \Z} Q_t^{(p)}$. 
We mention that, for $ \a\in (1,2)$ and $q\le 2(\a-1)/(2-\a)$, the two conditions on $(r_n)$ may possibly be different.
\par
It remains to derive the moment properties of $W^{(p)}=\sum_{t\in \Z} Q_t^{(p)}$, $p>\a$.  Fix an integer $m\ge 2$. 
The following identities may be of the type 
$\infty=\infty$ due to the lack of integrability. 
By the change of measure \eqref{eq:tildeQ} we have
\beao
\E\big[|W^{(p)}|^m]
&=&\E [ \|Q\|_p^{\a- m} \, \|Q\|_1^m  ] /\E\big[\|Q\|_p^\a\big]\,,
\eeao
where we used the fact that $\Theta_t\ge 0$ a.s. By the definition of $Q=\Theta/\|\Theta\|_\a$ it remains to 
verify whether the quantity
\beao
\E \big[ \|Q\|_p^{\a- m} \, \|Q\|_1^m\big]=
\E\big[\|\Theta\|_p^{\alpha-m}\,\|\Theta\|_1^m/\|\Theta\|_\alpha^\alpha\big]
\eeao
is finite or not.
This formula involves the backward tail process. Applying the time-change formula
for \regvary\ stationary processes (\cite{basrak:segers:2009} or
Proposition 4.3.1 in \cite{mikosch:wintenberger:2024}) we may conclude for $h\ge 0$
that
\beao
\P((\Theta_{-h},\ldots,\Theta_0)\in\cdot,\Theta_{-h}\neq 0)&=&\dfrac{\E[\sum_{j=h}^{K-1}\vep_{(\1(t\le j))_{0\le t\le h}}(\cdot)]}{\E[K]}\\
&=&\P((\Theta_{0},\ldots,\Theta_h)\in \cdot,\Theta_{h}\neq 0)\,.
\eeao
Therefore the backward and forward tail processes have the same \ds ,
$\Theta_t\in\{0,1\}$, $t\in \Z$, and 
\beao
\E\big[\|\Theta\|_p^{\alpha-m}\,\|\Theta\|_1^m/\|\Theta\|_\alpha^\alpha\big]&=&\E\Big[\Big(\sum_{t\in\Z}\Theta_t\Big)^{m-(m-\a)/p-1}\Big]\,.
\eeao
If $m-(m-\a)/p-1<0$ the \rhs\  is finite since $\sum_{t\in\Z}\Theta_t\ge 1$. This can happen only if $p\le 1$ because it implies $\a<p\le 1$ and then, choosing $m$ sufficiently large, one obtains the existence of every moment for $W^{(p)}$ and the limit ratio. 
\par
Therefore we focus on the case $p>1$ and $m-(m-\a)/p-1>0$. We have
\beao
\E\Big[\Big(\sum_{t\in\Z}\Theta_t\Big)^{m-(m-\a)/p-1}\Big]&\ge &\E\Big[\Big(\sum_{t\ge 0}\Theta_t\Big)^{m-(m-\a)/p-1}\Big]\\
&= &\dfrac{\E\Big[\sum_{j=0}^{K-1}\Big(\sum_{t\ge 0}\1(t\le j)\Big)^{m-(m-\a)/p-1}\Big]}{\E[K]}\\
&= &\dfrac{\E\Big[\sum_{j=1}^{K}j^{m-(m-\a)/p-1}\Big]}{\E[K]}\,.
\eeao
The right-hand side is infinite if and only if $m-(m-\a)/p\ge q$ 
  since we have required $\E[K^q]=\infty$.
As soon as $m-(m-\a)/p>(2(\a-1)/(2-\a))\vee 2$ it is possible to find such a $q$ and we rewrite the condition on $m$ as follows, 
\beao 
m>\dfrac{2(\a-1)/(2-\a)\vee 2-\a/p}{1-1/p}\,,
\eeao
thus providing a simplified expression for $\a\le 1.5$. 
{\red   Applying
Lemma \ref{lem:x} and \ref{lem:2}  for $\a\in(0,1)$ and $\a\in (1,2)$, respectively, we conclude that $R_\a^{(p)}$ has infinite moments of integer-valued or even order $m$, respectively. Recalling the definition of $R_\a^{(p)}$, we have  for $\a\in (1,2)$:
\beao
R_\a^{(p)}&:=&\dfrac{ \sum_{i=1}^\infty \big( \Gamma_i^{-1/\alpha} W_i^{(p)}-
C_i^{(\alpha)}\, \E[W^{(p)}]\big) }{\big\|\big(\Gamma_i^{-1/\alpha}\big)_{i\ge 1}\big\|_\infty }\\
&=&R_{\a,p}\dfrac{ \big\|\big(\Gamma_i^{-1/\alpha}\big)_{i\ge 1}\big\|_p}{\big\|\big(\Gamma_i^{-1/\alpha}\big)_{i\ge 1}\big\|_\infty }=:R_{\a,p}\wt R_p\,.
\eeao
For $\a\in(0,1)$ 
the relation $R_\a^{(p)}=R_{\a,p}\wt R_p$ remains valid since $W_i^{(p)}>0$ a.s. for every $i\ge1$.}
H\"older's inequality yields for any $r,s>1$ satisfying $1/r+1/s=1$,
\beao
\E[|R_\a^{(p)}|^m]\le \big(\E[|R_{\a,p}|^{mr}]\big)^{1/r}\,\big(\E\big[\textcolor{red}{\wt R_p^{ms}}\big]\big)^{1/s}\,.
\eeao
Considering $r=m'/m$, the result will follow if $\wt R_p$ has finite moments of any order. 
Using Corollary 4.4 of \cite{matsui:mikosch:wintenberger:2024}, we compute the Laplace transform
\textcolor{red}{
\beao
\E\big[\ex^{-\lambda \wt R_p}\big]=
\dfrac{\ex^{-\lambda}}{1+\int_0^{1}(1-e^{-y^p\lambda})d(-y^{-\a})}
=\dfrac{\ex^{-\lambda}}{1+\lambda^{\a/p}\int_0^{\lambda^{1/p}}(1-e^{-y^p})d(-y^{-\a})}\,,\qquad \lambda>0\,,
\eeao
which also exists on $\lambda\in [-\vep,0]$ for sufficiently small $o<\vep<1$.
Indeed, the denominator tends to $1$ as $\lambda \to 0$ and due to the inequality 
$c_1|x|<|e^x-1|<c_2|x|$ for $0<|x|<1$, we have for small $\lambda>0$, that  
\[
 \big[\ex^{\lambda \wt R_p}\big] \le 
\dfrac{\ex^{\lambda}}{1+c \lambda \int_0^{1}y^p d(-y^{-\a})}= \dfrac{e^{\lambda}}{1-c\lambda\a/(p-a)}<\infty, 
\]
where we keep the denominator positive.} Thus a similar reasoning as in the proof of Lemma \ref{lem:1} 
shows the finiteness of the moment generating function of $\wt R_p$ around its origin. 
Therefore $\wt R_p$ has finite moments of any order, and the desired result follows.
\section{Explicit formulas for even moments of $W^{(2)}$} \setcounter{equation}{0}\label{sec:even:moments}
In this section  we assume that $m$ is an even number and $(X_t)$ 
is a univariate \regvary\ stationary \seq\ with index \textcolor{red}{$\a\in (1,2)$}. 
We will use the \pp\ \rep s in \eqref{eq:series:rep} to derive the existence
and \rep\ of the $m$th moment of  $R_{\a,2}:=\xi_\a/\zeta_{\a,2}$.
Recall the definition of $W^{(2)}$ from \eqref{eq:suppress}.
\bpr\label{prop:m-moment}
Assume that $\alpha\in (1,2)$ and 
\beam\label{condi:mom:q}
\E \big[(W^{(2)})^m\big] <\infty\quad \mbox{ for some $m\in \{2,4,\ldots\}$.}
\eeam
Then  $\E[R_{\alpha,2}^m]$ is finite and has \rep
\beam\label{eq:simple2}
\dfrac{(m-1)!}{(m/2-1)!} {\sum}'
(k-1)! \dfrac 2\a
\Big(\dfrac{ \alpha/2}{\Gamma(1-\a/2)}\Big)^{k}
\prod_{j=1}^k \dfrac {\Gamma((\ell_j+1-\alpha)/2)} {\ell_j! \big(1\vee \big(m-j-\sum_{i=1}^j \ell_i \big)\big)}
\prod_{j=1}^{k}
\E\big[( W^{(2)})^{\ell_j+1}\big]\,, \nonumber\\
\eeam
where ${\sum}'$ denotes summation over the index set
\beao
\Big\{1\le k\le m\,, 0\le \ell_1,\ldots,\ell_k\le m: \ell_1+\cdots +\ell_k=m-k\Big\}\,.
\eeao
\epr
\bre
If \eqref{condi:mom:q} is not satisfied then the \rep\ \eqref{eq:simple2} fails. Indeed, for $k=1$ the product $\prod_{j=1}^{k}
\E\big[(W^{(2)})^{\ell_j+1}\big]$ would degenerate to $\E\big[(W^{(2)})^{m}\big]=\infty$. 
\ere
\bre For $m=2$, \eqref{eq:simple2} yields the formula
\beam\label{second:moment:R_a^2}
\E[R_{\alpha,2}^2] 
&=&  \E[(W^{(2)})^2] +(\E[ W^{(2)}])^2\,\dfrac{\alpha}{2}\, \Big(\dfrac{\Gamma((1-\alpha)/2)}{\Gamma(1-\alpha/2)}\Big)^2\,.
\eeam
From Matsui et al. \cite{matsui:mikosch:wintenberger:2024} we also know
\beao
\E[R_{\a,2}] =\E [W^{(2)}] \,\dfrac{\Gamma((1-\a)/2)}{\sqrt{\pi}\,\Gamma(1-\a/2)}\,.
\eeao
Hence 
\beam\label{eq:variance}
\var(R_{\a,2})=\E[(W^{(2)})^2] +(\E[W^{(2)}])^2\,\Big(\dfrac{\Gamma((1-\alpha)/2)}{\Gamma(1-\alpha/2)}\Big)^2
\Big(\dfrac{\a}2 -\dfrac 1{\pi}\Big)\,.
\eeam
\ere
 
%

\begin{proof}Recall from Lemma~\ref{prop:Wtilde} that for $p=2>\a$, $u\in\bbr$, $\la>0$, 
\beao
\Psi_X(u,\la)&:=&\E\big[\exp\big(i\,u\,\xi_\alpha-\lambda \,\zeta_{\alpha,2}^2\big)\big]\\
&=& \exp\Big(\E[\|Q\|_2^\alpha]\, 
\int_0^\infty \E\Big[\exp\Big(i\,y\,u  W^{(2)} -\lambda y^2\Big)-1-i\,y\,
u \,W^{(2)} 
\Big] d(-y^{-\alpha})\Big)\,.\nonumber
\eeao
Since we deal with the ratio $R_{\a,2}=\xi_\a/\zeta_{\a,2}$ the quantity
$\E[\|Q\|_2^\alpha]$ will not appear, and therefore we assume without loss of generality that it equals 1. Then we have
\beao
\dfrac{\partial \Psi_X(u,\lambda)}{\partial u}&=& 
\E\big[i\xi_\alpha \exp\big(i\,u\,\xi_\alpha-\lambda \,\zeta_{\alpha,2}^2\big)\big]=: \varphi_X(u,\lambda)\, \Psi_X(u,\lambda)\,,
\eeao
where 
\beao
\varphi_X(u,\lambda):= i \,\int_0^\infty \E\big[y\,  W^{(2)} \, 
\textcolor{red}{\big( \exp\big(i\,y\,u  W^{(2)} -\la \,y^2 \big)-1\big)}\big]\, d(-y^{-\alpha})\,. 
\eeao
Derivatives of $\varphi_X$ and $\Psi_X$ \wrt\ 
$u$ will be denoted by $\varphi^{(j)}$ and $\Psi_X^{(j)},\,j=1,2,\ldots$. 
\par
We need the following auxiliary result whose proof will be given later.
\ble\label{lem:m-moment}
Under the assumptions of Proposition~\ref{prop:m-moment}, 
$\E[\xi_\alpha^m \exp(-\lambda\, \zeta_{\alpha,2}^2)]<\infty$ for  every $\lambda>0$.
\ele
\noindent
Applying Leibniz's rule 
to $\varphi_X(u,\lambda) \,\Psi_X(u,\lambda)$, we obtain 
\beam\label{eq:Leibniz}
 \E\big[\xi_\alpha^m \,\exp\big(-\lambda \,\zeta_{\alpha,2}^2\big)\big]
&=&(-i)^m \,\Psi_X^{(m)}(0,\la)\nonumber\\
&=& (-i)^m \,\sum_{j=0}^{m-1} {{m-1}\choose j}\,
\varphi_X^{(j)}(0,\la)\, \Psi_X^{(m-1-j)}(0,\la)\,,
\eeam
where, by change of variable,  for every $1\le j <m$, 
\beao
 \varphi^{(j)}(0,\lambda) &=&i^{j+1}\, \int_0^\infty \E
\big[(y\, W^{(2)})^{j+1} \ex^{-\la\,y^2}\big]\, d(-y^{-\alpha}) \\
&=& i^{j+1} \E \big[
( W^{(2)})^{j+1} 
\big] \int_0^\infty y^{j+1} \ex^{-\la \,y^2} d(-y^{-\alpha})\\
&=& i^{j+1} \E \big[
( W^{(2)})^{j+1} 
\big] \int_0^\infty z^{(j+1)/2} \ex^{-\la \,z} d(-z^{-\alpha/2}) .
\eeao
{\red The following formula is derived from the definition of the $\Gamma$-function; it will be used frequently:
\beam\label{eq:gammaformula}
\dfrac1{\lambda^{1/p}}=\dfrac  1{\Gamma(1/p)}\int_0^\infty y^{1/p-1}\ex^{- \la y}dy=\dfrac  p{\Gamma(1/p)}\int_0^\infty\ex^{-\la y^p }dy\,,\qquad \la>0\,.
\eeam
Applying \eqref{eq:gammaformula} we obtain
\beao
 \varphi^{(j)}(0,\lambda) &= &i^{j+1}\, \E \big[
\big( W^{(2)}\big)^{j+1} 
\big]\, \dfrac{\alpha}{2}\, \lambda^{-(j+1-\alpha)/2}\, \Gamma((j+1-\alpha)/{2})\,.
\eeao}
Here we also used the moment condition 
$\E[( W^{(2)})^m]<\infty$; see 
\eqref{condi:mom:q}. This formula remains valid for $j=0$ although 
its derivation is slightly different. 
Thus \eqref{eq:Leibniz} turns into
\beao\lefteqn{(-i)^m \,\Psi_X^{(m)}(0,\la)}\\
&=&(-i)^m 
\sum_{j=0}^{m-1} {{m-1}\choose j } i^{j+1}\, \E \big[
\big(W^{(2)} \big)^{j+1}\big] \,
\dfrac{\alpha}{2} \,\la^{-(j+1-\alpha)/{2}} 
\Gamma((j+1-\alpha)/{2}) \,\Psi_X^{(m-1-j)}(0,\lambda)\,. 
\eeao
Iterating this step for $\Psi_X^{(\cdot)}(0,\lambda)$, we arrive at 
the following expression, where ${\sum}'$ stands for summation over the 
indices introduced in the proposition:
\beao
&&(-i)^m \,\Psi_X^{(m)}(0,\la)\\
&=& {\sum}'
 {{m-1}\choose {\ell_1}} \, {{m-2-\ell_1} \choose {\ell_2}} \cdots  
{{m-k-\sum_{i=1}^{k-1}\ell_i} \choose {\ell_k}} \\
&&\quad  \times \E \big[
\big( W^{(2)}\big)^{\ell_1+1} \big]\cdots \E \big[
\big( W^{(2)}\big)^{\ell_k+1} \big]\,(\alpha/{2})^k\, \la^{-(\ell_1+1-\a)/{2}-\cdots-(\ell_k+1-\a)/{2}}\,\\
&&\quad \times \Gamma((\ell_1+1-\alpha)/{2})\cdots \Gamma((\ell_k+1-\alpha)/{2})\, 
\Psi_X(0,\lambda) \\
&=& {\sum}'
\dfrac{(m-1)! }{\prod_{j=1}^{k} \big(\ell_j !\,
(1\vee (m-j-\sum_{i=1}^j \ell_i ))\big)}  
\prod_{j=1}^k \E\big[( W^{(2)} )^{\ell_j+1} \big]\,(\alpha/{2})^k \prod_{j=1}^k \Gamma((\ell_j+1-\alpha)/{2})\\ &&
\times\underbrace{\la^{(k\alpha-m)/{2}}\,\Psi_X(0,\lambda)}_{=:g(\la)}\,. 
\eeao
From {\red Formula \eqref{eq:gammaformula}} and Fubini's theorem we may conclude that 
\beam
\E\big[R_{\a,2}^m\big]&=&\E\Big[\dfrac{\xi_\a^m}{\zeta_{\a,2}^m}\Big]\\
&=&\dfrac{1}{\Gamma(m/{2})}\E\Big[\xi_\alpha^m \int_0^\infty \la^{m/2-1}\,\ex^{-\lambda \zeta_{\alpha,2}^2}\,d\la\Big]\nonumber\\
&=&\dfrac{1}{\Gamma(m/{2})}\int_0^\infty \la^{m/2-1} 
\E\big[\xi_\alpha^m \,\ex^{-\lambda \zeta_{\alpha,2}^2}\big]\,d\la
\,.\label{eq:finale}
\eeam
{\red Using again Formula \eqref{eq:gammaformula}} we also obtain
\beao
 \dfrac{1}{\Gamma(m/2)} \int_0^\infty g(\la)\,\la^{m/2-1}\,d\la
&=&\dfrac{1}{\Gamma({m}/{2})} \int_0^\infty
\lambda^{k\alpha/{2}-1}\, \Psi_X(0,\lambda)\,d\lambda\\ 
&=&  \dfrac{1}{\Gamma({m}/{2})} \int_0^\infty \lambda^{k\alpha/{2}-1} 
\exp\big(-\la^{\alpha/2}\, \E\|Q\|_2^{\alpha}\,\Gamma(1-\alpha/2)\big)\,d\lambda \\
&=&\dfrac{2}{\alpha} \dfrac{\Gamma(k)}{\Gamma(m/2)} 
\dfrac{1}{\big(\E\|Q\|_2^{\alpha} \Gamma(1-\alpha/2)\big)^k} \,.
\eeao
Combining this calculation with \eqref{eq:finale}, \eqref{eq:Leibniz} and recalling 
that we have assumed without loss of generality that  $\E\|Q\|_2^{\alpha}=1$,
we have proved the desired relation \eqref{eq:simple2}.
\end{proof}

\begin{proof}[Proof of Lemma~\ref{lem:m-moment}]
We consider the 
series representation of $\xi_\alpha$ given in \eqref{measure:change:series:expression} for $p=2$ where we assume without loss of generality that the constant 
$\E[\|Q\|_2^\a]=1$.
Then we have the following  decomposition of $\xi_\a$: 
\beao
 \xi_\alpha &=& \sum_{i>m/\alpha} \Gamma_i^{-1/\alpha} ( W^{(2)}_i-\E[W^{(2)}])
+  \E[W^{(2)}]\,
\sum_{i>m/\alpha} (\Gamma_i^{-1/\alpha}-C_i^{(\alpha)}) \\
&& + \sum_{i\le m/\alpha} \Gamma_i^{-1/\alpha}  W^{(2)}_i
- \sum_{i\le m/\alpha} C_i^{(\alpha)} 
\E[ W^{(2)}]
=: I_1+ \textcolor{red}{ \E[ W^{(2)}]\,I_2+ I_3}-C_4,  
\eeao
where the infinite series $I_1,I_2$ are finite a.s.; see Samorodnitsky and 
Taqqu~\cite[p.30]{samorodnitsky:taqqu:1994}). 
Since $m$ is an even integer it suffices to show that $\E[I_1^m+I_2^m+I_3^m\,\ex^{-\la \zeta_{\a,2}^2}]<\infty$, $\la>0$.
\par
We start by bounding $\E [I_2^m]$. By a generalized H\"older inequality
we obtain for some constant $c$,
\beam
 \E [I_2^m] &\le& \E\Big[\Big(
\sum_{i>m/\alpha} (\Gamma_i^{-1/\alpha}-C_i^{(\alpha)})
\Big)^m\Big]\nonumber \\
&=& \sum_{i_1,\ldots,i_m >m/\alpha} \E\Big[
\prod_{j=1}^m \Big|\Gamma_{i_j}^{-1/\alpha}-C_{i_j}^{(\alpha)}\Big|
\Big]\nonumber\\
&\le& c \sum_{i_1,\ldots,i_m >m/\alpha} 
\prod_{j=1}^m \Big( \E \big[\big(\Gamma_{i_j}^{-1/\alpha}-C_{i_j}^{(\alpha)}\big)^m \big] 
\Big)^{1/m}. \label{eq:ineq:binom}
\eeam
A binomial expansion 
yields for $i>m/\a$,
\beam\label{eq:binom:exp}
 \E \Big[\big(\Gamma_i^{-1/\alpha}-C_i^{(\alpha)}\big)^m\Big] &=& 
\sum_{\ell =0}^m {{m}\choose {\ell}}\,\E \big[\Gamma_i^{-\ell/\alpha}\big] 
(-C_i^{(\alpha)})^{m-\ell} \nonumber\\
& =& \sum_{\ell =0}^m {{m}\choose {\ell}} 
\dfrac{\Gamma(i-\ell/\alpha)}{\Gamma(i)} (-C_i^{(\alpha)})^{m-\ell}\,. 
\eeam
By Stirling's formula for the gamma function as $i\to\infty$,
\beao
 \dfrac{\Gamma(i-\ell/\alpha)}{\Gamma(i)} &\sim& \sqrt{\dfrac{i}{i-\ell/\alpha}} 
\dfrac{(i-\ell/\alpha)^{i-\ell/\alpha}}{i^i} \ex^{\ell/\alpha}\, (
1+O(i^{-1})) \\
&=& \ex^{\ell/\alpha} \,\Big(1-\dfrac{\ell/\alpha}{i}\Big)^i\, (i-\ell/\alpha)^{-\ell/\alpha}\,(
1+O(i^{-1})) \\
&=& (
1+O(i^{-1})) \,i^{-\ell/\alpha}\,. 
\eeao
In the last step we used an expansion for real $x$ as $i\to\infty$, 
\beao
 \Big(
1+\dfrac{x}{i}
\Big)^i = \ex^x\, \Big(
1-\dfrac{x^2}{2i}+ \dfrac{x^3 (8+11\,x)}{24 \,i^2}+\cdots\Big)\,. 
\eeao
Thus we have uniform bounds for $\ell=1,\ldots,m$, small $\vep>0$ and sufficiently large $i$,
\beao
 (1-\vep)\, i^{-\ell/\alpha} \le \dfrac{\Gamma(i-\ell/\alpha)}{\Gamma(i)} \le  
(1+\vep)\, i^{-\ell/\alpha}\,. 
\eeao
\textcolor{red}{Inserting the lower (upper) bounds to the negative (positive) terms 
in \eqref{eq:binom:exp}, we conclude that
\beao
 \E \big[\big(\Gamma_i^{-1/\alpha}-C_i^{(\alpha)}\big)^m\big] \le (1+\vep) \big(i^{-1/\alpha}-C_i^{(\alpha)}\big)^m \,,\qquad i>m/\a\,.
\eeao}
\textcolor{blue}{This part seems to be wrong and insufficient.}
In view of \eqref{eq:ineq:binom} we get for some constant $c$,
\beao
 \E [I_2^m] \le c \,\sum_{i_1,\ldots,i_m > m/\alpha} \prod_{j=1}^m |i_j^{-1/\alpha}-C_{i_j}^{(\alpha)}| \le c \,\Big(
\sum_{i>m/\alpha} \big|i^{-1/\alpha}-C_i^{(\alpha)}\big|
\Big)^m <\infty\,. 
\eeao
\par
Next we will bound 
\beao
\E\big[I_3^m \ex^{-\lambda \zeta_{\alpha,2}^2}\big] = \E\Big[\Big(
\sum_{i\le m/\alpha} \Gamma_i^{-1/\alpha}  W_i^{(2)}\,\exp\Big(-(\lambda/m)\, \sum_{k\ge 1} \Gamma_k^{-2/\alpha} \Big)
\Big)^m\Big]\,.
\eeao
Since $c_0:=\sup_{y>0}y\,\exp(-(\la/m) y^2)<\infty$ we have for every $i\ge 1$, 
\beao
\Gamma_{i}^{-1/\alpha}   \exp\Big(-(\lambda/m) 
\sum_{k\ge 1} \Gamma_k^{-2/\alpha} \Big)<
 \Gamma_{i}^{-1/\alpha} \exp\big(-(\lambda/m) \Gamma_i^{-2/\alpha}\big)\le c_0\quad \as
\eeao
Hence for some constant $c$,
\beao
 \E\big[I_3^m \,\ex^{-\lambda \,\zeta_{\alpha,2}^2}\big] 
&\le& c_0^m \E \Big[\Big(\sum_{i\le m/\a} \big| W_i^{(2)}\big|\Big)^m \Big]\\
&\le & c\, (m/\a)^{m-1}\,\E\big[\big|W^{(2)}\big|^m\big]\,.
\eeao
The \rhs\ is finite in view of our assumption \eqref{condi:mom:q}. 
\par
It remains to bound $\E [I_1^m]$. Since $(W_i^{(2)}-\E[W^{(2)}])$ is
iid centered and $f(x)=x^m$ is convex for even $m$ we can use the symmetrization 
Lemma~6.3 in Ledoux and Talagrand \cite{ledoux:talagrand:1991}:
for a Rademacher \seq\ $(\vep_i)$ independent of $(W_i^{(2)})$ and $(\Gamma_i)$ we have for some 
constants $c$, possibly different from line to line,
\beao
\E\big[I_1^m\big]&=&
\E\Big[\E\Big[\Big(\sum_{i>m/\a} \Gamma_i^{-1/\a}\,( W_i^{(2)}-\E[ W^{(2)}])\Big)^m\,\Big|\,(\Gamma_i)\Big]\Big]\\
&\le & c\,\E\Big[\E\Big[\Big(\sum_{i>m/\a} \Gamma_i^{-1/\a}\,( W_i^{(2)}-\E[ W^{(2)}])\,\vep_i\Big)^m\,\Big|\,(\Gamma_i)\Big]\Big]\\
&\le &c\,\E\Big[\max_{i> m/\a} (i/\Gamma_i)^{m/\a}\Big]\,
\E\Big[\Big(\sum_{i>m/\a} i^{-1/\a}\,\vep_i\,(W_i^{(2)}-\E[W_i^{(2)}])\Big)^m\Big].
\eeao
In the last step we used the contraction principle for  a Rademacher series (Theorem~4.4 in  
Ledoux and Talagrand \cite{ledoux:talagrand:1991}) conditional on $(\Gamma_i)$.
The second expectation on the \rhs\ is finite since the summands in the corresponding infinite series are independent and have $m$th momemt (Theorem 6.11 in  \cite{ledoux:talagrand:1991}). For $\E[I_1^m]<\infty$ it remains to show that $V=\E[\max_{i> m/\a} (i/\Gamma_i)^{m/\a}]$ is finite.
We have for any constant $c>e$ that
\beao
V\le \E\big[\max_{i>m/\a} (i/\Gamma_i)^{m/\a}\1(i/\Gamma_i\le c) \big]+ \E\big[\max_{i>m/\a} (i/\Gamma_i)^{m/\a}\1(i/\Gamma_i> c) \big]=:V_1+V_2\,.
\eeao
We have $V_1\le c^{m/\a}$ while by the choice $c>e$ we obtain for some constants $c_1,c_2, c_3$,
\beao
V_2&\le& \sum_{i>m/\a} i^{m/\a}\,\E [\Gamma_i^{-m/\a}\1(i/\Gamma_i> c)]
= \sum_{i>m/\a} i^{m/\a}\,\int_0^{i/c} x^{-m/\a+i-1}\,\ex^{-x}\,dx/\Gamma(i)\\
&\le &c_1\,\sum_{i>m/\a} i^{m/\a}\,(i/c)^{i-m/\a}/\Gamma(i)
\le c_2\,\sum_{i>m/\a} i^{m/\a}\,(i/c)^{i-m/\a}/ ((i/e)^{i}(2\pi\,i)^{1/2})\\
& \le & c_3\sum_{i>m/\a} (e/c)^{i} (2\pi\,i)^{-1/2}<\infty\,.
\eeao
In the last steps we exploited a series \rep\ of the incomplete gamma \fct\ (Gradshteyn and Ryzhik \cite{gradshteyn:ryzhik:1980}, p. 941) and Stirling's formula for the gamma \fct .
Thus we have proved the lemma.

\end{proof}
\bexam \rm Consider a \regvary\ linear causal process 
$(X_t)$ with index $\a\in (1,2)$; see Section \ref{sec:linear}. Then we have the explicit expressions
\beao
Q^{(2)}&=& \dfrac{\Theta}{\|\Theta\|_2}=\Theta_Z\,\Big(\dfrac{\psi_{t+J}}{\|\psi\|_2}
\Big)_{t\in\bbz}\,,\\
W^{(2)}&=&\sum_{j\in\bbz}Q^{(2)}= \Theta_Z\,\dfrac{\sum_{j\in\bbz} \psi_j}{\|\psi\|_2}
=:\Theta_Z\,c_2\,,\qquad 
R_{\a,2}^X:=\dfrac{\xi_\a^X}{\zeta_{\a,2}^X}= c_2\,\dfrac{\xi_\a^Z}{\zeta_{\a,2}^Z}=:c_2\,R_{\a,2}^Z\,, 
\eeao
\eqref{eq:simple2} simplifies: 
\beao
\E\big[(R_{\a,2}^X)^m\big]=c_2^m\,\E\big[(R_{\a,2}^Z)^m\big]\,,
\eeao
and $\E\big[(R_{\a,2}^Z)^m\big]$ is obtained from \eqref{eq:simple2} by replacing $\prod_{j=1}^k \E\big[(W^{(2)})^{\ell_{j+1}}\big]$ with
\beam\label{eq:qm}
\prod_{1\le j\le k\,,\ell_1+\cdots+\ell_k=m-k}
\E[(\Theta_Z)^{\ell_j+1}]
&=&\prod_{1\le j\le k\,,\ell_1+\cdots+\ell_k=m-k} \big(q_++ q_-\, (-1)^{\ell_{j}+1}\big)\nonumber\\
&=& (q_+-q_-)^{n_k}\,,
\eeam
where $n_k=n_k(\ell_1,\ldots,\ell_k)$ is the number of even $\ell_j$ among $\ell_1,\ldots,\ell_k$ \st\  $\ell_1+\cdots +\ell_k=m-k$ and, if none of
them is odd, $(q_+-q_-)^{0}$ is interpreted as 1. In the extreme cases when $q_+=1$ or $q_-=1$, \eqref{eq:qm} equals 1 or $(-1)^{n_k}=(-1)^m=1$ for even $m$, respectively.
If $q_+=q_-=0.5$,   \eqref{eq:qm} equals 0 or 1 according as $n_k>0$ or $n_k=0$.
\par
We observe that $W^{(2)}$ degenerates to $\Theta_Z$ for the $Z$-\seq .
In view of \eqref{eq:variance} we get an expression for the variance of $R_{\a,2}^X$:
\beao
\var(R_{\a,2}^X)=c_2^2 \var(R_{\a,2}^Z)=c_2^2\,\Big(1 + (q_+-q_-)^2\,
\Big(\dfrac{\Gamma((1-\alpha)/2)}{\Gamma(1-\alpha/2)}\Big)^2
\Big(\dfrac{\a}2 -\dfrac 1{\pi}\Big)
\Big)\,.
\eeao
If $q_+=q_-=0.5$ the \rhs\ degenerates to $c_2^2$ and does not depend on $\a$. It is not difficult to see that this is the \asy\ variance
for $(S_n/\gamma_{n,2})$ in the case of a Gaussian linear process with the same coefficients $(\psi_j)$.
\par 
In particular,
for a \regvary\  AR(1) process $(X_t)$ with parameter 
$\varphi\in (-1,1)$ we obtain
\beao
\E\big[(R_{\a,2}^X)^m\big]=\Big(\dfrac{\sqrt{1+\varphi}}{\sqrt{1-\varphi}}\Big)^m\,\E\big[(R_{\a,2}^Z)^m\big]\,.
\eeao

\par
To grasp further intuition on the formula \eqref{eq:simple2}
we illustrate how we can calculate the even $m$th 
moments in the iid symmetric case. Then, with our previous notation,
$c_2=1$ and $q_\pm=1/2$. 
In Table \ref{table:1} for each $m$, we present non-zero combinations of $(k,(\ell_i)_{i\le k})$ in the sum 
\eqref{eq:simple2} together with the values of $\E[(R_{\alpha,2}^Z)^m],\,m=2,4,6,8$. 
\begin{center}
\begin{tabular}{|c|c|c|c|}\hline
$m$ &$k$ & $(\ell_1,\ldots,\ell_k)$ & $\E[(R_{\alpha,2}^Z)^m]$ \rule[-8pt]{0pt}{20pt} \\ \hline 
$2$  & $1$ & $\ell_1=1$ & $1$  \\ \hline
\multirow{2}{*}{4}   & $2$ & $\ell_1=\ell_2=1$ & \multirow{2}{*}{$1+\alpha$} \\
    & $1$ & $\ell_1=3$& \\ \hline
\multirow{4}{*}{6}  & $3$   & $\ell_1=\ell_2=\ell_3=1$ &\multirow{3}{*}{$1+3\alpha+2\alpha^2$ }\\
   & $2$   & $(\ell_1,\ell_2)=(1,3),(3,1)$ &\\
   & $1$   & $\ell_1=5$ &\\ \hline
\multirow{4}{*}{8}    & $4$   & $\ell_1=\ell_2=\ell_3=\ell_4=1$ &\multirow{4}{*}{$\frac{1}{3}(3+20\alpha+34\alpha^2+17\alpha^3)$} \\
    & $3$   & $(\ell_1,\ell_2,\ell_3)=(1,1,3),(1,3,1),(3,1,1)$ &\\
    & $2$  & $(\ell_1,\ell_2)=(3,3),(1,5),(5,1)$ &\\
    & $1$   & $\ell_1=7$ & \\ \hline
\end{tabular}\vspace{1mm}
\captionof{table}{$\E[(R_{\alpha,2}^Z)^m],\,m=2,4,6,8$.}\label{table:1}
\end{center}
Our results coincide with those of Logan et al. \cite[(5.21)]{logan:mallows:rice:shepp:1973}. In the latter paper the following equation is basic
\beam\label{eq:bothsides}
 \frac{1}{\pi} \int_0^\infty \E[\ex^{iR_{2,\alpha}t}] \ex^{-st} dt = \int_0^\infty \ex^{-s^2 t^2/2} \mathscr{D}(t) dt,\,s>0, 
\eeam
where $\mathscr{D}(t)$ is defined in \eqref{eq:parab};
see Appendix~\ref{calc:density:R:iid} below for further explanations. 
To obtain expressions for the moments 
\cite[(5.21)]{logan:mallows:rice:shepp:1973} used Taylor expansions on the both sides of \eqref{eq:bothsides} followed by a comparison of the  
coefficients.
\eexam
\begin{comment}
\end{document}